\documentclass[oneside,a4paper,11pt,notitlepage]{article}
\usepackage[left=0.8in, right=0.8in, top=1.5in, bottom=1.5in]{geometry}
\usepackage[T1]{fontenc}
\usepackage[utf8]{inputenc}
\usepackage{tocloft}

\AtBeginDocument{\selectfont}

\usepackage[english]{babel}
\usepackage{float}
\usepackage{amssymb}
\usepackage{lmodern}
\usepackage{amsmath}
\usepackage{amsthm}
\usepackage{bm}
\usepackage{mathtools}
\usepackage{braket}
\usepackage{esint}
\usepackage{enumerate}
\usepackage{booktabs}
\usepackage{graphicx}
\usepackage{tikz}
\usetikzlibrary{patterns}
\usepackage{multicol}
\usepackage{caption}
\usepackage{comment}
\usepackage[hyphens]{url}
\usepackage{xcolor}
\usepackage{aliascnt}
\usepackage[maxbibnames=99]{biblatex}
\usepackage{hyperref}
\usepackage{csquotes}

\hypersetup{
  colorlinks=true,
  linkcolor={magenta},
  citecolor={cyan},
  urlcolor={cyan},
  breaklinks=true,
  bookmarksnumbered,
  pdftitle={A common quantitative approach to Chiti and Kohler--Jobin inequalities}
}

\allowdisplaybreaks[2]
\theoremstyle{plain}
\newtheorem{theorem}{Theorem}
\numberwithin{theorem}{section}
\numberwithin{equation}{section}

\newaliascnt{proposition}{theorem}
\newtheorem{proposition}[proposition]{Proposition}
\aliascntresetthe{proposition}

\newaliascnt{lemma}{theorem}
\newtheorem{lemma}[lemma]{Lemma}
\aliascntresetthe{lemma}

\newaliascnt{corollary}{theorem}
\newtheorem{corollary}[corollary]{Corollary}
\aliascntresetthe{corollary}

\newaliascnt{definition}{theorem}
\theoremstyle{definition}
\newtheorem{definition}[definition]{Definition}
\aliascntresetthe{definition}

\newaliascnt{example}{theorem}

\aliascntresetthe{example}

\newaliascnt{remark}{theorem}
\theoremstyle{definition}
\newtheorem{remark}[remark]{Remark}
\aliascntresetthe{remark}

\addto\extrasenglish{%
}

\newcommand{\R}{\mathbb R}

\newcommand{\abs}[1]{\left|#1\right|}
\newcommand{\norm}[1]{\left\lVert#1\right\rVert}
\newcommand{\Tmod}{T_{\mathrm{mod}}}
\newcommand{\res}{\mathcal R}

\newcommand{\Addresses}{{ 
 \bigskip

  \noindent\textit{E-mail address}, A.L.~Masiello: \texttt{albalia.masiello@unina.it} \par\nopagebreak \par\nopagebreak 
   \medskip

 \noindent\textit{E-mail address}, C.~Nitsch: \texttt{c.nitsch@unina.it} \par\nopagebreak

  \medskip 
 
 \noindent\textit{E-mail address}, C.~Trombetti : \texttt{cristina@unina.it} \par\nopagebreak 
   \medskip

 \noindent\textsc{Dipartimento di Matematica e Applicazioni ``R. Caccioppoli'', Universit\`a degli studi di Napoli Federico II, Via Cintia, Complesso Universitario Monte S. Angelo, 80126 Napoli, Italy.}\par\nopagebreak 
 
 \bigskip
 
 \noindent\textit{E-mail address}, V. ~Amato: \texttt{v.amato@ssmeridionale.it} \par\nopagebreak 

     \medskip 

      \noindent\textit{E-mail address}, L. ~Barbato: \texttt{l.barbato@ssmeridionale.it} \par\nopagebreak 

     \medskip 

 \noindent \textit{E-mail address}, F.~Villone: \texttt{f.villone@ssmeridionale.it} \par\nopagebreak 
  
   \medskip 
 
  \noindent\textsc{Mathematical and Physical Sciences for Advanced Materials and Technologies, Scuola Superiore Meridionale, Largo San Marcellino 10, 80138 Napoli, Italy. } \par\nopagebreak 
  \bigskip

 \textit{E-mail address}, E.~Cristoforoni: \texttt{emanuele.cristoforoni@unimi.it} \par\nopagebreak 
  \medskip 
\noindent\textsc {Dipartimento di Matematica ``Federigo Enriques'', Universit\'a degli Studi di Milano La Statale, Via Saldini 50 20123 Milano Italia.}\par\nopagebreak

}}

\title{A unified quantitative approach to Chiti and Kohler--Jobin inequalities}
\author{V. Amato, L. Barbato, E. Cristoforoni, A.L. Masiello, C. Nitsch, \\
C. Trombetti, F. Villone}
\date{}

\begin{document}
\maketitle

\begin{abstract}
We obtain two quantitative results related to classical rearrangement inequalities for the first Dirichlet eigenvalue. The first one is a quantitative version of Chiti's reverse H\"older inequality. More precisely, for every $1\le q<p<\infty$, we prove that the deficit controls, up to translations, the $L^q$ distance between a first eigenfunction of a bounded domain and the first eigenfunction of the ball having the same first eigenvalue. The second result concerns the Kohler--Jobin inequality and gives a quantitative estimate in terms of the $L^\infty\!\to L^2$ distance between the corresponding Dirichlet resolvents. In both cases, the main idea is to find a suitable superlevel set of the first eigenfunction which is almost optimal for the Faber--Krahn inequality.


\end{abstract}

\hypersetup{
                linkcolor={black},
}

\begin{center}
\begin{minipage}{10cm}
\small
\tableofcontents
\end{minipage}
\end{center}

\hypersetup{
                linkcolor={magenta},
}

\section{Introduction}

Symmetrization techniques have been extensively used to obtain qualitative information on eigenvalues and eigenfunctions of elliptic operators. Among the classical results in this direction, we recall the two famous isoperimetric results for the eigenvalue and the torsional rigidity of the Dirichlet Laplacian, the Faber-Krahn inequality and the Saint-Venant inequality. If $\Omega$ is an open bounded set in $\R^n$, we recall that the first Dirichlet eigenvalue $\lambda_1(\Omega)$ is the smallest value for which the problem

\begin{equation}
    \label{eigenvalue:intro}
    \begin{cases}
        \displaystyle-\Delta u = \lambda_1(\Omega) u& \text{in}\;\Omega\\
u = 0 & \text{on}\;\partial\Omega, 
    \end{cases}
\end{equation}
admits a non-trivial solution, while the Torsional rigidity $T(\Omega)$ is the $L^1$ norm of the solution to  

\begin{equation*}
    \begin{cases}
        \displaystyle-\Delta v = 1& \text{in}\;\Omega\\
v = 0 & \text{on}\;\partial\Omega. 
    \end{cases}
\end{equation*}

The Schwarz symmetrization is the common tools that allows to prove that among sets of fixed volume, the first Dirichlet eigenvalue $\lambda_1(\Omega)$ is minimized, while the Torsional rigidity is maximized by the ball.\medskip 

Of course, these are just two of the many famous examples of isoperimetric properties proved by means of symmetrization techniques. 
The reverse H\"older inequalities for eigenfunctions of the Dirichlet Laplacian are other examples. Let $u$ be a first eigenfunction of \eqref{eigenvalue:intro}, then  it is possible to prove that for all $p> q$

\begin{equation}
    \label{chiti:intro}
    \int_\Omega u^p\, dx\le K(n,\lambda_1(\Omega), p, q) \left(\int_\Omega u^q\, dx\right)^\frac{p}{q}, 
\end{equation}
where the constant $K(n, \lambda_1(\Omega), p, q)$ is explicitly given by
\begin{equation*}
   K(n,\lambda_1(\Omega), p,q)=\dfrac{\displaystyle{\int_{B^\lambda} z^p\, dx}}{\displaystyle{\left(\int_{B^\lambda} z^q\, dx\right)^\frac{p}{q}}},
\end{equation*}
where $B^\lambda$ is the ball such that $\lambda_1(B^\lambda)=\lambda_1(\Omega)$ and $z$ is a first eigenfunction of $B^\lambda$. This result was firstly proved by Payne and Rayner \cite{PR72} in dimension $n=2$  and then generalized in any dimension by Kohler--Jobin \cite{Kohler81} and Chiti \cite{ChitiIso} in higher dimension with different techniques. 
In particular, it is worth mentioning that reverse H\"older inequalities hold more generally for eigenfunctions of uniformly elliptic second-order operators with homogeneous Dirichlet boundary conditions;   see \cite{ChitiRH, ChitiIso}.
These results have been also extended in several directions. For instance, Alvino, Ferone and Trombetti \cite{AFT98} considered nonlinear eigenvalue problems involving the $p$-Laplacian. A reverse H\"older inequality for extremal Sobolev functions was proved in \cite{CarrollRatzkin15}. More recently, Chiti-type inequalities and their stability have also been studied in the setting of $\mathrm{RCD}(K,N)$ spaces, see \cite{GunesMondino23,Wu24,wu25}.\medskip 

When it comes to isoperimetric results proved by means of rearrangement techniques, one cannot avoid mentioning the celebrated Kohler--Jobin inequality, which settles the P\'olya--Szeg\H{o} conjecture:

\begin{quote} \it
   among all sets with prescribed torsional rigidity, the ball minimizes the principal frequency.
\end{quote}

In closed form, the inequality proved in \cite{Kohler78} can be summarized as follows

\begin{equation}
    \label{KJ:intro}
    T(\Omega)^\frac{2}{n+2}\lambda_1(\Omega)\ge  T(B)^\frac{2}{n+2}\lambda_1(B),
\end{equation}
and it is proved by using a symmetrization that preserves the Dirichlet integral and increases the $L^p$ norm of a function (see also \cite{KohlerJobin82}). A modern presentation and several extensions of this method can be found in \cite{Brasco14}. Further developments include the reverse Kohler--Jobin problem \cite{BrianiButtazzoGLB24}, the Gaussian case \cite{HerscoviciLivshyts24}, and recent $p$-Kohler--Jobin inequalities in Euclidean and Gaussian spaces \cite{CFMW26,CFMWGauss26}. We also mention the related study of functionals involving the first eigenvalue and torsional rigidity in \cite{AmatoNitschTrombettiVillone26}.\medskip

 Chiti's reverse H\"older inequality and the Kohler--Jobin inequality belong to the same classical rearrangement framework and their proofs use closely related level-set techniques. 
In Chiti's argument, the concentration of the first eigenfunction $u$ on $\Omega$ is compared with that of the first radial eigenfunction $z$ on the ball $B^\lambda$. The Kohler--Jobin proof also starts from the first eigenfunction, but constructs a radial rearrangement preserving the Dirichlet integral, using the torsional rigidity to calibrate the comparison.\medskip 
  
The aim of the present paper is twofold: on one side, we highlight the similarities between the two problems; on the other side, we prove strengthened versions of inequalities \eqref{chiti:intro} and \eqref{KJ:intro}, in the spirit of the recent literature on the stability of symmetrization inequalities.

The pioneering result for Schwarz symmetrization is the quantitative P\'olya-Szeg\H{o} principle by Cianchi, Esposito, Fusco and Trombetti \cite{CEFT08}, which controls the distance of a Sobolev function from its Schwarz symmetrization in terms of the deficit in the Dirichlet integral. In \cite{ABMP24}, it was shown that this result can be used to derive stability estimates for comparison principles obtained via rearrangement, and in particular to prove a quantitative version of Talenti's comparison result \cite{T76}. This strategy turned out to be quite flexible and has been adopted in several other contexts: comparison results for first-order Hamilton-Jacobi equations \cite{AB25}, for the $k$-Hessian equation \cite{MS25}, and for Talenti-type inequalities with Robin boundary conditions \cite{ABCMP26}, as well as other problems in which radial symmetry plays a key role, such as \cite{CL257, DGPV25}.\medskip

In all these works, the problems under consideration are studied under a measure constraint or, more generally, a geometric constraint. Schwarz rearrangement plays a crucial role in the proof of these results, and its quantitative counterpart, namely the quantitative P\'olya--Szeg\H{o} principle \cite{CEFT08}, is the key ingredient for obtaining stability estimates.\medskip

The situation for the Chiti and the Kohler--Jobin inequalities is slightly different: even if we are still considering optimization problems for which a ball is the solution, for \eqref{chiti:intro} and \eqref{KJ:intro} we are imposing a spectral constraint rather than a geometric one, so this has to be taken into account if one wants to prove a stability result. By a "stability result" we mean an inequality of the form

\begin{equation*}
   D_{\mathrm{Ch}}(\Omega):=\dfrac{\displaystyle{\int_{B^\lambda} z^p\, dx}}{\displaystyle{\left(\int_{B^\lambda} z^q\, dx\right)^\frac{p}{q}}}-\dfrac{\displaystyle{\int_{\Omega} u^p\, dx}}{\displaystyle{\left(\int_\Omega u^q\, dx\right)^\frac{p}{q}}}\ge C\, d(\Omega, B^\lambda)
\end{equation*}
for the Chiti inequality, and
\begin{equation*}
     \delta_{\mathrm{KJ}}(\Omega):=T(\Omega)^{\frac{2}{n+2}}\lambda_1(\Omega)-T(B)^{\frac{2}{n+2}}\lambda_1(B)\ge C\, d(\Omega, B^T)
\end{equation*}
for the Kohler--Jobin one, where $d$ is a suitable distance between sets.\medskip

In the existing literature, stability is typically measured through the Fraenkel asymmetry
\begin{equation*}
    \mathcal{A}(\Omega):=\min_{x\in\mathbb{R}^n}\left\{\frac{|\Omega\,\Delta\, B_r(x)|}{|\Omega|}\, :\, |B_r(x)|=|\Omega|\right\}.
\end{equation*}
Somehow surprisingly, the Fraenkel asymmetry is not the right quantity for a quantitative version of these two inequalities. Indeed, in \autoref{sec:counterexamples}, we exhibit two sequences of sets for which, respectively,  $D_{\mathrm{Ch}}$ and $\delta_{\mathrm{KJ}}$ tend to zero, while the Fraenkel asymmetry remains bounded away from zero.\medskip

The necessity to overcome this issue leads us to change approach. We would like to mention that a source of inspiration was the work by G\"une\c{s} and Mondino \cite{GunesMondino23}, where they prove a reverse H\"older inequality for eigenfunctions on $\mathrm{RCD}(K,N)$-spaces, together with a quantitative issue which, however, is different from ours. What is clear from their argument is that, in order to obtain a quantitative Chiti inequality, it is crucial to show that a set $\Omega$ for which the Chiti inequality is almost an equality can be replaced, at a small cost, by a suitable subset, $U$, of $\Omega$ having the same measure as the reference ball $B^\lambda$. 
Once we define $U$, we can focus our attention on a quantitative Faber-Krahn inequality. The failure of the Fraenkel asymmetry leads us to consider the \emph{spectral stability} of the Faber-Krahn inequality  proved by Allen, Kriventsov and Neumayer in \cite{AKN23, AKN25,AKN21} (see also \cite{NTV2026} for the Robin case), instead of the classical, geometric form of the inequality, proved in \cite{BDPV15} (see \autoref{AKN-res} and \autoref{thm:eigenfunction:fk}).\medskip

In particular \autoref{thm:eigenfunction:fk} states that the Faber--Krahn deficit controls the $L^2$-distance of the normalized eigenfunction of $\Omega$ from that of the ball of the same volume. This version will be used in the quantitative Chiti inequality, which we will state as follows. 
Let $B^\lambda=B^{\lambda_1(\Omega)}$ be the ball such that
\[
 \lambda_1(B^\lambda)=\lambda_1(\Omega),
\]
and let $u$ and $z$ be positive first eigenfunctions of $\Omega$ and $B^\lambda$, respectively, extended by zero outside their domains.
We define the scaling-invariant Chiti deficit
\begin{equation}
\label{delta:chiti}
     \delta_{\mathrm{Ch}}(\Omega):
 ={\lambda_1(\Omega)^\frac{n(q-p)}{2q}}\left( \frac{\|z\|_{L^p(B^\lambda)}^p}{\|z\|_{L^q(B^\lambda)}^p}-\frac{\|u\|_{L^p(\Omega)}^p}{\|u\|_{L^q(\Omega)}^p}\right)\ge0.
\end{equation}

We can now state our first main theorem.

\begin{theorem}\label{thm:intro-Chiti}
Let $1\le q<p<\infty$. There exists a positive constant $C=C(n,p,q)$ such that, for every bounded, open, and connected set $\Omega\subset\mathbb R^n$ and every non-negative first eigenfunctions $u$ and $z$  on $\Omega$ and $B^\lambda$ respectively,
\[
    \delta_{\mathrm{Ch}}(\Omega)\geq C\left(\inf_{x_0\in\mathbb R^n}{\norm{\frac{u}{\norm{u}_{q}}-\frac{z(\,\cdot-x_0)}{\norm{z}_{q}}}_{q}^{\max\{2,q\}}}\right)^{p+1}.
\]
\end{theorem}

The proof follows the classical comparison argument of Chiti. We then consider the level $\bar t$ such that the measure of the superlevel set $U=\{u>\bar t\}$ is exactly $\abs{B^\lambda}$, and we prove that the Chiti deficit forces $\bar t$ to be small. Combining this information with the techniques developed in the quantitative Talenti problem \cite{ABMP24}, we show that $U$ has a small Faber-Krahn deficit. At this point, the eigenfunction stability estimate of \cite{AKN23} gives the desired conclusion.\medskip

\autoref{AKN-res}, on the other hand, states that the Faber--Krahn deficit controls the distance in operator norm of the resolvent operator on $\Omega$, $\res_\Omega$,  from that of a ball of the same volume (see \autoref{spectraFK} for the exact definition of the resolvent operator). This version of the resolvent stability enters the proof of the quantitative Kohler--Jobin inequality.
Let $B^T=B^{T(\Omega)}$ be the ball such that $T(B^T)=T(\Omega)$. We prove the following result.

\begin{theorem}\label{thm:intro-KJ}
Let $n\ge2$. There exists a positive, dimensional constant $c_n$ such that, for every bounded, open, and connected set $\Omega\subset\mathbb R^n$,
\[
 \delta_{\mathrm{KJ}}(\Omega)\ge c_n \left(\inf_{x\in\R^n}\dfrac{\norm{\res_\Omega-\res_{B^{T(\Omega)}(x)}}_{L^\infty\to L^2}}{T(\Omega)^\frac{n+4}{2(n+2)}}\right)^4.
\]
\end{theorem}

 Again, to prove \autoref{thm:intro-KJ}, we start from the classical proof of \cite{Brasco14} and we identify a suitable subset of $\Omega$  that is almost optimal for Faber--Krahn inequality, so that the resolvent stability estimate of \cite{AKN25} gives the conclusion.\medskip

Although the two quantitative proofs follow the same line, the two stability statements are substantially different, and this difference is intrinsic to the nature of the inequalities. Chiti's inequality involves different norms of the first eigenfunction and therefore, essentially, \emph{only} detects the first spectral mode. On the other hand, in the Kohler--Jobin inequality, the presence of the Torsional rigidity makes the quantitative inequality in \autoref{thm:intro-KJ} sensitive to \emph{all} of the spectrum of $-\Delta$.  This is linked to the notion of \emph{$\gamma$-convergence} (see \autoref{defi:discussion-gamma}). Indeed,  \autoref{thm:intro-KJ} can also be stated in the following way

\begin{theorem}\label{thm-KJ-gamma}
Let $n\ge2$. There exists a positive, dimensional constant $C_n$ such that, for every bounded, open set $\Omega\subset\mathbb R^n$,
\begin{equation}\label{eq:discussion-KJ-gamma-quantitative:intro}
\delta_{\mathrm{KJ}}(\Omega)\ge C_n\left(\inf_{x\in\mathbb R^n}
\frac{d_\gamma(\Omega,B^{T(\Omega)}(x))}{T(\Omega)}\right)^4.
\end{equation}

Consequently, if $(\Omega_k)$ is a sequence of bounded open sets
and
\[
r_k=\left(\frac{T(\Omega_k)}{T(B_1)}\right)^{1/(n+2)},
\]
then
\[
\delta_{\mathrm{KJ}}(\Omega_k)\longrightarrow0
\]
if and only if there exist translations $x_k\in\mathbb R^n$ such that
\[
\frac{\Omega_k-x_k}{r_k}
\xrightarrow{\gamma}B_1.
\]
\end{theorem}

The relation between the two statements \autoref{thm:intro-KJ} and \autoref{thm-KJ-gamma} is discussed at the end of \autoref{KJ-fc}. 
We remark that the connectedness assumption in \autoref{thm:intro-Chiti} and \autoref{thm:intro-KJ} is only needed to simplify the discussion. The cases of a disconnected set is discussed in \autoref{cor:Chiti-disconnected} and \autoref{cor:KJ-disconnected}.\medskip

The paper is organized as follows. In Section 2, we introduce the notation and the preliminary results that will be used throughout the paper. Section 3 is devoted to the quantitative Chiti inequality in the $L^q$ setting, while in Section 4 we prove the quantitative Kohler--Jobin inequality. Finally, in Section 5 we discuss examples illustrating the different rigidity mechanisms of the two inequalities and the limitations of geometric stability without additional assumptions.

\medskip\noindent\textbf{Disclosure of AI assistance.}
The authors used generative AI tools during the preparation of this manuscript.
For the quantitative Chiti inequality, AI was used only for writing, editing, and improving the presentation of arguments developed and verified by the authors.
For the quantitative Kohler–Jobin inequality, the authors used ChatGPT 5.6 Sol High also at an exploratory stage, asking it to adapt the strategy of the quantitative Chiti proof to the Kohler–Jobin functional.
Since the two inequalities are closely related, the model was instructed to follow the same proof structure, with particular attention to \autoref{lem:Chiti-core-eigenvalue}, where the two arguments differ, and this led to the preliminary formulation of \autoref{lem:KJ-exact} and \autoref{prop:tbar-KJ}. The prompt and the raw AI-generated output are available at: 

\begin{center}
    \href{https://github.com/vincamato/quantitative-chiti-Kohler--Jobin}{https://github.com/vincamato/quantitative-chiti-Kohler--Jobin}
\end{center}

The resulting preliminary manuscript was then carefully reviewed, corrected, and mathematically verified by the authors.

AI tools were also used to assist with the writing and editing of the paper.
The authors take full responsibility for the mathematical content of the final manuscript.

\section{Notation and preliminary results}

Throughout this article, $|E|$ denotes the Lebesgue measure of a measurable set $E\subset\mathbb R^n$. If the space is clear from the context, we will simply write $\|\cdot\|_r$ for the $L^r$ norm. All the functions defined on a domain will be extended by zero outside the domain whenever necessary.

\medskip
We start by fixing the notation for eigenvalues, torsional rigidity and their scaling.
Let $D\subset\mathbb R^n$ be a bounded open set. We denote by
\begin{equation*}
    \lambda_1(D):=\min_{\varphi\in H_0^1(D)\setminus\set{0}}\dfrac{\displaystyle\int_D|\nabla\varphi|^2\,dx}{\displaystyle\int_D\varphi^2\,dx}\quad\text{and}\quad T(D):=\max_{\varphi\in H_0^1(D)\setminus\set{0}}\dfrac{\displaystyle\left(\int_D\varphi\,dx\right)^2}{\displaystyle\int_D|\nabla\varphi|^2\,dx}.
\end{equation*}
If $w_D$ is the torsion function, $-\Delta w_D=1$ in $D$ and $w_D\in H_0^1(D)$, then
\[
 T(D)=\int_D w_D\,dx=\int_D|\nabla w_D|^2\,dx.
\]
For $r>0$, we have
\[
 \lambda_1(rD)=r^{-2}\lambda_1(D),
 \qquad
 T(rD)=r^{n+2}T(D).
\]
Thus the quantity $\lambda_1(D)T(D)^{2/(n+2)}$ is scale invariant.\medskip

We recall that the Faber-Krahn and Saint-Venant inequalities state that, for every bounded open set $D\subset\mathbb{R}^n$ and every ball $B\subset\mathbb{R}^n$,
\begin{equation*}
    |D|^{\frac{2}{n}}\lambda_1(D)\ \ge\ |B|^{\frac{2}{n}}\lambda_1(B),
\end{equation*}
\begin{equation*}
    T(D){|D|^{-\frac{n+2}{n}}}\ \le T(B){|B|^{-\frac{n+2}{n}}}.
\end{equation*}
The Kohler--Jobin inequality reads
\begin{equation*}
 T(D)^\frac{2}{n+2}\lambda_1(D)
 \ge T(B)^\frac{2}{n+2}\lambda_1(B),
\end{equation*}
with equality only for balls, up to the usual negligible-set convention.\medskip

Finally we also recall the Hong--Krahn--Szeg\H{o} inequality (see \cite{P55}).  Let $B$ be a ball, then for every bounded open set $E$ with $\abs{E}=\abs{B}$, we have
\begin{equation}\label{eq:HKS}
 \lambda_2(E)\ge 2^{\frac2{n}}\lambda_1(B).
\end{equation}
Here and in what follows, eigenvalues are always counted with multiplicity.

 \subsection{Properties of rearrangements}
Let us recall some basic notions about rearrangements. For a general overview, see, for instance, \cite{kes}.

 \begin{definition}\label{distribution:function}
	Let $f: \Omega \to \R$ be a measurable function, the \emph{distribution function} of $f$ is the function $\mu : [0,+\infty[\, \to [0, +\infty[$ defined as the measure of the superlevel sets of $\abs{f}$, i.e.,
	$$
	\mu(t)= \abs{\Set{x \in \Omega \, :\,  \abs{f(x)} > t}}.
	$$
\end{definition}
\begin{definition}\label{decreasing:rear}
	Let $f: \Omega \to \R$ be a measurable function, the \emph{decreasing rearrangement} of $f$, denoted by $f^\ast(\cdot)$, is defined as
$$f^*(s)=\inf\{t\geq 0:\mu(t)<s\}.$$
\end{definition}
	\begin{remark}
	Let us notice that the function $\mu(\cdot)$ is decreasing and right continuous, and the function $f^\ast(\cdot)$ is its generalized inverse.
    From Definitions \ref{distribution:function} and \ref{decreasing:rear}, one can prove that
	 $$f^\ast (\mu(t)) \leq t, \quad \forall t\ge 0,$$ 
  $$\mu (f^\ast(s)) \leq s \quad \forall s \ge 0.$$

We will mostly apply the decreasing rearrangement to first eigenfunctions $u$ of a bounded, \emph{connected} open set $\Omega\subset\mathbb{R}^n$. In that case,  the function $u$ is analytic in $\Omega$ and cannot be constant on any open set, which forces $\mu$ to be continuous and $f^\ast$ to be its inverse.

	\end{remark}

	\begin{definition}
	 The \emph{Schwarz rearrangement} of $f$ is the function $f^\sharp $ whose superlevel sets are concentric balls with the same measure as the superlevel sets of $f$. 
	\end{definition}

		We have the following relation between $f^\sharp$ and $f^*$:
	$$f^\sharp (x)= f^*(\omega_n\abs{x}^n),$$
 and one can easily check that the functions $f$, $f^*$, and $f^\sharp$ are equi-distributed, i.e., they have the same distribution function, and there holds
\begin{equation*}
\displaystyle{\norm{f}_{L^p(\Omega)}=\norm{f^*}_{L^p(0, \abs{\Omega})}=\lVert{f^\sharp}\rVert_{L^p(\Omega^\sharp)}}, \quad \text{for all } p\ge1.    
\end{equation*}
Finally, we recall the contraction property
\begin{equation*}
 \norm{f^\sharp-g^\sharp}_{L^r(\R^n)}
 \le \norm{f-g}_{L^r(\R^n)},
 \qquad 1\le r<\infty,
\end{equation*}
for non-negative functions extended by zero.

\subsection{Chiti comparison theorem}

The following consequence of Chiti's comparison theorem will be useful in the sequel \cite{ChitiRH, ChitiIso}.

\begin{theorem}[Chiti, \cite{ChitiRH}]\label{comparison:chiti}
Let $\Omega\subset\R^n$ be a bounded, open set and let $B^\lambda$ be the ball satisfying $\lambda_1(B^\lambda)=\lambda_1(\Omega)$.  Let $u,z>0$ be first Dirichlet eigenfunctions on $\Omega$ and $B^\lambda$ respectively, normalized by \begin{equation}\label{eq:normalization:q}
\norm{u}_{q}=\norm{z}_{q},  \text{ for some } q\ge1
\end{equation}
  Extend $z^*$ by zero to $(|B^\lambda|,|\Omega|)$.  Then
\begin{equation}\label{eq:Chiti-majorization}
 F(s):=\int_0^s\big((z^*)^q-(u^*)^q\big)\,dt\ge0
 \qquad(0\le s\le|\Omega|).
\end{equation}

In particular, for $p>q\ge1$ one has 
\begin{equation}
\label{norme:p}
\norm{u}_{p}\le\norm{z}_{p}.    
\end{equation}
\end{theorem}

\begin{proof}[Sketch of the proof]
Since the map $\tau\mapsto\tau^{p/q}$ is convex and of class $C^1$ in $[0,+\infty)$, for every $t\in[0,|\Omega|]$ we have
\[
    (z^*)^p-(u^*)^p\ \ge\ \frac pq\,(u^*)^{p-q}\big((z^*)^q-(u^*)^q\big)=\frac pq\,(u^*)^{p-q}F'(t),
\]
with
\[
    F(t)=\int_0^t (z^*)^q-(u^*)^q\,dt .
\]
Note that $F(0)=F(|\Omega|)=0$ by the normalization \eqref{eq:normalization:q}. Integrating by parts,
\[
    \int_0^{|\Omega|}\frac pq\,(u^*)^{p-q}F'(t)\,dt=\int_0^{|\Omega|}\frac pq\,\Big(-(u^*)^{p-q}\Big)'F(t)\,dt ,
\]
and $\big(-(u^*)^{p-q}\big)'\ge0$ because $u^*$ is decreasing and $p>q$. Hence, proving \eqref{eq:Chiti-majorization} will also imply \eqref{norme:p}.

We are left to show that $F(t)\ge0$ for every $t\in[0,|\Omega|]$. 

To do so, 
we first show that
\begin{equation}\label{claim1}
    \norm{u}_{L^\infty(\Omega)}\le\norm{z}_{L^\infty(B^\lambda)},
    \qquad\text{i.e.}\qquad u^*(0)\le z^*(0).
\end{equation}
so that $F'(t)\ge 0$ in a right neighborhood of $0$. Then
we show that there exists a unique $s^*$ such that
\begin{equation}
\label{claim2}
   F'(s^*)=\big((z^*)^q-(u^*)^q\big)(s^*)=0 . 
\end{equation} 
Certainly, there exists a point $s_1$ at which it vanishes. The fact that this point is unique will give \eqref{eq:Chiti-majorization}.To prove both \eqref{claim1} and \eqref{claim2}, we argue by contradiction.\medskip

Suppose by contradiction that $u^*(0)>z^*(0)$. By continuity, $u^*>z^*$ in a right neighbourhood of $0$; let
\[
    S:=\sup\{s\in(0,|B^\lambda|]:\ u^*>z^* \text{ in } [0,s)\}.
\]
If $S=|B^\lambda|$, then $u^*>z^*$ in $[0,|B^\lambda|)$ and $u^*\ge0=z^*$ in $(|B^\lambda|,|\Omega|)$, so that
\[
    \int_0^{|\Omega|}(u^*)^q\,dt>\int_0^{|\Omega|}(z^*)^q\,dt,
\]
which contradicts \eqref{eq:normalization:q}. If $S<|B^\lambda|$, then $u^*(S)=z^*(S)$ and
\begin{equation}\label{eq:crossing0}
    u^*(s)>z^*(s)\qquad\forall s\in[0,S).
\end{equation}
Define
\[
    w(s)=\begin{cases}
        u^*(s), & s\in[0,S),\\
        z^*(s), & s\in[S,|B^\lambda|].
    \end{cases}
\]
By Talenti's theorem \cite{T76}, $w$ satisfies
\[
    -w'(s)\le \lambda_1(\Omega)\,\frac{s^{-2+\frac2n}}{n^2\omega_n^{\frac2n}}\int_0^s w(t)\,dt ,
\]
and the inequality is strict in $[S,|B^\lambda|)$ by \eqref{eq:crossing0}.

Hence, it can be shown that the function $x\mapsto w(\omega_n|x|^n)$ belongs to $W^{1,2}_0(B^\lambda)$ and has Rayleigh quotient strictly smaller than $\lambda_1(\Omega)=\lambda_1(B^\lambda)$, which is absurd. This proves \eqref{claim1}.

To prove \eqref{claim2}, we suppose by contradiction that there exists another point $s_2\in(0,|B^\lambda|)$ $s_2\neq s_1$ such that $u^*(s_2)=z^*(s_2)$ and
\begin{equation*}
    u^*(s)>z^*(s)\qquad\forall s\in(s_1,s_2).
\end{equation*}
Define
\[
    w(s)=\begin{cases}
        z^*(s), & s\in[0,s_1]\cup[s_2,|B^\lambda|],\\
        u^*(s), & s\in(s_1,s_2).
    \end{cases}
\]

Therefore $w$ satisfies
\[
    -w'(s)\le \lambda_1(\Omega)\,\frac{s^{-2+\frac2n}}{n^2\omega_n^{\frac2n}}\int_0^s w(t)\,dt .
\]

As before, the function $x\mapsto w(\omega_n|x|^n)$ belongs to $W^{1,2}_0(B^\lambda)$ and satisfies
\[
    \dfrac{\displaystyle\int_{B^\lambda}|\nabla w|^2\,dx}{\displaystyle\int_{B^\lambda}w^2\,dx}<\lambda_1(\Omega),
\]
which is absurd.
\end{proof}

\autoref{comparison:chiti}  is also known as mass concentration comparison (see, for instance, \cite{ALT89}). Let us stress, that in particular, the following comparison holds
\begin{equation}\label{eq:Chiti-uniform-bound}
 \norm{u}_\infty\le\norm{z}_\infty,
\end{equation}
under the usual normalization \eqref{eq:normalization:q}.

\subsection{The Kohler--Jobin rearrangement}

One of the main tools for the Kohler--Jobin rearrangement is the modified torsional rigidity with respect to a reference function. We refer the reader to \cite{Brasco14,KohlerJobin82} for an in-depth discussion on the topic.

\begin{definition}
    In the notation of \cite{Brasco14}, we say that a non-negative function $w\in H^1_0(\Omega)\cap L^\infty(\Omega)$ is a reference function if
\[t\mapsto\dfrac{\abs{\set{w>t}}}{\displaystyle\int_{\set{w=t}} \abs{\nabla w}\,d\mathcal H^{n-1}}\in L^\infty(\,[0,\norm{w}_\infty]).\]
Then, if $w$ is a reference function, we define the modified torsional rigidity as 
\begin{equation*}
 \Tmod(\Omega;w)
 :=\sup\Set{
 \dfrac{\left(\displaystyle\int_\Omega g(w)\,dx\right)^2}{\displaystyle\int_\Omega\abs{\nabla g(w)}^2\,dx}| g\in\operatorname{Lip}([0,\norm w_\infty]),\ g(0)=0,\ g(w)\neq0}. 
\end{equation*}
\end{definition}
Notice that
\begin{equation*}
 \Tmod(\Omega;w)\le T(\Omega)
\end{equation*}
for every reference function $w$.\medskip

Let $M$ denote the $L^\infty$-norm of the eigenfunction $u$ and, for every $t\in(0,M)$ let 
\[\Omega_t=\set{u>t},\qquad\mu(t)=\abs{\Omega_t}.\]
Then, by \cite[Lemma~3.2 and Proposition~3.5]{Brasco14} we have the following proposition
\begin{proposition}
    For every $t\in(0,M)$ the function $(u-t)_+$ is a reference function on $\Omega_t$, moreover letting
    \[\mathcal T(t) = \Tmod(\Omega_t,(u-t)_+),\]
    we have
    \begin{equation}\label{eq:Tt-formula}
 \mathcal T(t)=\int_t^M\frac{\mu(s)^2}{\lambda_1(\Omega)F(s)}\,ds,
\end{equation}
where
\[F(t)=\int_{\Omega_t} u\,dx.\]
\end{proposition}
We can now introduce the Kohler--Jobin rearrangement of the eigenfunction $u$. We remark that, actually, the Kohler--Jobin rearrangement is well defined for every reference function $w$ (see \cite[Lemmata 2 and 3]{KohlerJobin82} and \cite[Proposition 4.1]{Brasco14}).  
\begin{theorem}\label{KJ-rearrangment}
    Let $B_r$ be the ball centred at the origin such that
    \[T(B_r)=\Tmod(\Omega,u).\]
    There exists a non-negative, radially symmetric decreasing function $v\in H^1_0(B_r)$ such that 
    \begin{equation}\label{eq:KJ-energy-equal}\int_{B_r} \abs{\nabla v}^2\,dx=\int_\Omega \abs{\nabla u}^2\,dx,\end{equation}
    and
    \begin{equation}\label{eq:KJ-L2-gain}\int_{B_r} v^2\,dx\ge \int_\Omega u^2\,dx.\end{equation}

\end{theorem}
\begin{remark}
    The previous theorem immediately yields the Kohler--Jobin inequality. Indeed, we have that
    \[\int_{B_r}\abs{\nabla v}^2\,dx=\int_\Omega \abs{\nabla u}^2\,dx=\lambda_1(\Omega),\qquad\int_{B_r}v^2\,dx\ge \int_\Omega u^2\,dx=1.\]
    Hence
    \[
        T(B_r)^{\frac{2}{n+2}}\,\lambda_1(B_r)\le T(B_r)^{\frac{2}{n+2}}\,\displaystyle\int_{B_r}\abs{\nabla v}^2\,dx = \Tmod(\Omega,u)^{\frac{2}{n+2}}\,\lambda_1(\Omega)\le T(\Omega)^{\frac{2}{n+2}}\, \lambda_1(\Omega).
\]
    
\end{remark}
\begin{proof}[sketch of the proof]
    The main idea of the proof is to rewrite the integral of the $L^2$-norm of $\nabla u$ in the variable $\mathcal{T} (t)$. That is, let $\varphi$ denote the inverse of the monotone function $\mathcal{T}$, then, using \eqref{eq:Tt-formula} and the coarea formula, after some computations, one can show that
    \[\int_\Omega \abs{\nabla u}^2\,dx= \int_0 ^{\mathcal{T}(0)} (\mu(\varphi(s))\,\varphi'(s))^2\,ds.\]
    On the other hand, given a decreasing function $\psi\,\colon[0,\mathcal{T}(0)]\to[0,+\infty)$, we can define the rearrangement $v$ to take the value $\psi(s)$ on the boundary of the ball whose torsional rigidity is exactly $s$, that is
    \[v(x)=\psi\left(\left(\frac{\abs{x}}{r_0}\right)^{n+2}\right).\]
    By direct computations, we then have 
    \[\int_{b_r} \abs{\nabla v}^2\,dx= \int_0 ^{\mathcal{T}(0)} \left(\abs{\set{v>\psi(s)}}\,\psi'(s)\right)^2\,ds,\]
    where, by construction, $\set{v>\psi(s)}$ is exactly the ball whose torsional rigidity is $s$ (and is hence explicit and independent of the choice of $\psi$). Then to have \eqref{eq:KJ-energy-equal} is enough to choose $\psi$ as the solution to
    \begin{equation*}
    \begin{cases}
        -\psi'(s)= -\varphi'(s)\, \dfrac{\mu(\varphi(s))}{\abs{\set{v>\psi(s)}}}\\[10pt]
        \psi(\mathcal{T}(0))=0.
    \end{cases}\end{equation*}
    Let $\Phi(t):=\psi\left(\varphi^{-1}(t)\right)$, then the torsional rigidity of the  ball $\set{v>\Phi(t)}$ is exactly $\mathcal{T}(t)$, hence, by the Saint--Venant inequality, \[\beta(t):=\abs{\set{v>\Phi(t)}}\le \mu(t).\]
    Hence, for every $t\in[0,M]$
    \[\Phi'(t)=\dfrac{\psi'\left(\varphi^{-1}(t) \right)}{\varphi'\left(\varphi^{-1}(t) \right)}=\dfrac{\mu(t)}{\beta(t)}\ge 1\]
    and, by the boundary condition $\Phi(0)=\psi(\mathcal{T}(0))=0$,
    \begin{equation}\label{Phigeid}\Phi(t)\ge t.\end{equation}
    Finally, from equations \eqref{Phigeid} and the identity $\beta(t)\Phi'(t)=\mu(t)$, one can prove \eqref{eq:KJ-L2-gain} using Cavalieri's principle and a change of variable, namely
    \[\begin{split}\int_{B_r} v^2\,dx&=2\int_0^{\norm{v}_\infty} \tau \abs{\set{v>\tau}}\,d\tau=2\int_0^{M}\Phi(t)\beta(t)\Phi'(t)\,dt\\[10pt]&=2\int_0^{M}\Phi(t)\mu(t)\,dt\ge 2\int_0^{M}t\mu(t)\,dt= \int_{\Omega} u^2\,dx.\qedhere \end{split}\]
\end{proof}
\subsection{Spectral stability of the Faber-Krahn inequality}\label{spectraFK}
We close the section with the two sharp quantitative results proved in \cite{AKN23,AKN25}. To state the result we need to introduce some terminologies.\medskip

For an open bounded set $\Omega\subset\mathbb{R}^n$ and a function $f\in H^{-1}(\mathbb{R}^n)$, let $u_\Omega^f\in H_0^1(\Omega)$ be the unique weak solution to the Poisson equation
\begin{equation*}
    \begin{cases}
        -\Delta u_\Omega^f = f & \text{in } \Omega,\\[5pt]
        u_\Omega^f = 0 & \text{on } \partial\Omega.
    \end{cases}
\end{equation*}
As for eigenfunctions, we extend $u_\Omega^f$ by zero to be defined as a function in $H_0^1(\mathbb{R}^n)$. The \emph{resolvent operator} is then defined as
\[
    \mathcal{R}_\Omega \colon f\in H^{-1}(\mathbb{R}^n)\mapsto u_\Omega^f\in H_0^1(\mathbb{R}^n).
\]

In \cite{AKN23} the authors prove the following quantitative estimate.
\begin{theorem}[ \cite{AKN25} Theorem 1.2 ]\label{AKN-res}
    Let $n\ge 2$. There exists a positive, dimensional constant $c_n$ such that for every bounded, open set $\Omega\subset\R^n$, with $\abs{\Omega}=\abs{B_1}$,
    \begin{equation}\label{eq:AKN-res}
    \lambda_1(\Omega)-\lambda_1(B_1) \ge\ c_n\norm{\mathcal{R}_\Omega-\mathcal{R}_{B_1(x_\Omega)}}_{L^\infty\to L^2}^2,
\end{equation}
for some point $x_\Omega\in\R^n$. In other words, for any $f\colon\mathbb{R}^n\to\mathbb{R}$ with $\norm{f}_{\infty}\le 1$, we have
\begin{equation*}
     \lambda_1(\Omega)-\lambda_1(B_1) \ge\ c_n\int_{\mathbb{R}^n}\abs{u_\Omega^f-u_{B_1(x_\Omega)}^f}^2\,dx .
\end{equation*}
\end{theorem}
We observe that \eqref{eq:AKN-res} implies that there exists a positive dimensional constant $C_n$ such that, for every bounded open set $\Omega\subset\R^n$,  
\[\abs{\Omega}^{\frac2{n}}\lambda_1(\Omega)-\abs{B}^{\frac2{n}}\lambda_1(B) \ge\ C_n\dfrac{\norm{\mathcal{R}_\Omega-\mathcal{R}_{B(x_\Omega)}}_{L^\infty\to L^2}^2}{\abs{\Omega}^\frac{n+4}{n}},\]
where $B$ is any ball such that $\abs{\Omega}=\abs{B}$.\medskip

From the previous theorem the authors deduce
\begin{theorem}[ \cite{AKN21} Theorem 1.1 ]\label{thm:eigenfunction:fk}
    Let $n\ge 2$. There exists a positive, dimensional constant $c_n$ such that for every bounded, open set $\Omega\subset\R^n$, with $\abs{\Omega}=\abs{B_1}$,
    \begin{equation}\label{eigenfunction:fk}
    \lambda_1(\Omega)-\lambda_1(B_1) \ge\ c_n\int_{\mathbb{R}^n}\abs{u_{\Omega}-z_{B(x_\Omega)}}^2\, dx,
\end{equation}
for some point $x_\Omega\in\R^n$. Here $u_\Omega$ and $z_B$ are two $ L^2$-normalized eigenfunctions, extended by zero to be defined on the whole $\mathbb{R}^n$.
\end{theorem}
We observe that \eqref{eigenfunction:fk} implies that there exists a positive dimensional constant $C_n$ such that, for every bounded open set $\Omega\subset\R^n$,  
\[\abs{\Omega}^{\frac2{n}}\lambda_1(\Omega)-\abs{B}^{\frac2{n}}\lambda_1(B) \ge C_n\int_{\mathbb{R}^n}\abs{u_{\Omega}-z_{B(x_\Omega)}}^2\, dx,\]
where $B$ is any ball such that $\abs{\Omega}=\abs{B}$ and  $u_\Omega$ and $z_B$ are two $ L^2$-normalized eigenfunctions.

\section{Quantitative Chiti inequality in \texorpdfstring{$L^q$}{Lq}}
In this section, we prove \autoref{thm:intro-Chiti}. Let $\Omega\subset\R^n$ be a bounded, open, and connected set, and fix $1\le q<p<\infty$. We write $\lambda=\lambda_1(\Omega)$ and denote by $B^\lambda$ the ball centred at the origin such that
\[
 \lambda_1(B^\lambda)=\lambda.
\]
For simplicity, throughout the proof we assume that
\begin{equation}\label{eq:Chiti-normalization}
 \lambda_1(\Omega)=1.
\end{equation}
The reduction to this normalization is justified at the end of the proof of \autoref{thm:intro-Chiti}.
Let $u$ and $z$ be the positive first eigenfunctions of $\Omega$ and $B^\lambda$, respectively, normalized by
\begin{equation}\label{eq:q-normalization}
 \norm{u}_{L^q(\Omega)}=\norm{z}_{L^q(B^\lambda)}=1.
\end{equation}
The scale-invariant Chiti deficit defined in \eqref{delta:chiti}, then reads
\begin{equation*}
 \delta:=\delta_{\mathrm{Ch}}(\Omega)=\norm{z}_{L^p(B^1)}^p-\norm{u}_{L^p(\Omega)}^p.
\end{equation*}
 We set $B:=B^1$, so that $B$ is the fixed ball with first eigenvalue equal to one.
All functions are extended by zero outside their domains. Unless otherwise specified, the constants in this section depend only on $n,p,q$.\medskip

The section is structured as follows. Following Chiti's comparison Theorem \ref{comparison:chiti}, we derive quantitative estimates for the decreasing rearrangements of $u$ and $z$. We then introduce the reference set
\[
 \bar t:=u^\ast(\abs{B}),\qquad
 U:=\{u>\bar t\},\qquad v:=(u-\bar t)_+,
\]
whose volume is $\abs{B}$, and prove that, for sufficiently small deficit,
\[
 \bar t\lesssim\delta^{1/(p+1)},\qquad
 0\le\lambda_1(U)-1\lesssim\delta^{1/(p+1)}.
\]
The quantitative Faber-Krahn inequality in the form \eqref{eigenfunction:fk} then gives
\[
 \norm{v-z(\,\cdot-x_\Omega)}_2^2
 \lesssim\delta^{1/(p+1)}
\]
for some $x_\Omega\in\R^n$. We conclude by passing from $L^2$ to $L^q$ on the supports of these functions and estimating the truncation error $u-v$.


\medskip
We begin with two elementary estimates. The first quantifies the comparison between the $p$th and $q$th powers; the second will allow us to pass from $L^2$ estimates to other finite norms on sets of uniformly bounded measure.

\begin{proposition}\label{lem:powers-Lq}
 For every $p>q\ge1$, there exists $c=c(p,q)>0$ such that, for all $a,b\ge0$,
 \begin{equation}\label{dis:pq}
  a^p-b^p
  \ge\frac pq\,b^{p-q}(a^q-b^q)
  +c(a-b)^2(a+b)^{p-2}.
 \end{equation}
 At $a=b=0$, the last term is defined by continuity.
\end{proposition}

\begin{proof}
 If $b=0$, the assertion holds for any $c\le1$. If $b>0$, setting $t=a/b$ reduces \eqref{dis:pq} to a positive lower bound for
 \[
  G(t):=\frac{t^p-1-\frac pq(t^q-1)}{(t-1)^2(t+1)^{p-2}},
  \qquad t\ne1.
 \]
 The strict convexity of $s\mapsto s^{p/q}$ shows that $G(t)>0$ for $t\ne1$. Moreover,
 \[
  G(0)=\frac pq-1>0,\qquad
  \lim_{t\to1}G(t)=\frac{p(p-q)}{2^{p-1}}>0,
  \qquad\lim_{t\to\infty}G(t)=1.
 \]
 Thus the continuous extension of $G$ to $[0,\infty)$ has a positive lower bound depending only on $p,q$.
\end{proof}

\begin{proposition}\label{lem:change-norm-Lq}
 Suppose that $h\in L^2(\R^n)\cap L^\infty(\R^n)$, $\abs{\{h\ne0\}}\le V$, and $\norm{h}_\infty\le K$. Then
 \[
  \begin{aligned}
   \norm{h}_r^2&\le V^{2/r-1}\norm{h}_2^2,
   &&1\le r\le2,\\
   \norm{h}_r^r&\le K^{r-2}\norm{h}_2^2,
   &&2\le r<\infty.
  \end{aligned}
 \]
\end{proposition}

\begin{proof}
 The first estimate follows from H\"older's inequality on the support of $h$. The second follows by integrating $\abs{h}^r\le K^{r-2}\abs{h}^2$.
\end{proof}

We can now estimate the distance between the decreasing rearrangements in terms of the Chiti deficit.

\begin{lemma}
 Under the normalizations \eqref{eq:Chiti-normalization} and \eqref{eq:q-normalization}, there exists $C=C(n,p,q)>0$ such that
 \begin{align}
  \int_0^{\abs{\Omega}}
  \abs{z^\ast-u^\ast}^2(z^\ast+u^\ast)^{p-2}\,ds
  &\le C\delta,\label{eq:weighted-Lq}\\
  \int_0^{\abs{B}}\abs{z^\ast-u^\ast}^2\,ds
  &\le C\delta^{\min\{1,2/p\}},\label{eq:core-L2-Lq}\\
  \int_{\abs{B}}^{\abs{\Omega}}(u^\ast)^q\,ds
  &\le C\delta^{\min\{1/2,q/p\}}.\label{lemma3.3last}
 \end{align}
\end{lemma}

\begin{proof}
 Applying \autoref{lem:powers-Lq} with $a=z^\ast(s)$ and $b=u^\ast(s)$ and integrating, we obtain
 \[
  \begin{split}
   \delta=\int_0^{\abs{\Omega}}\bigl((z^\ast)^p-(u^\ast)^p\bigr)\,ds\ge{}&c\int_0^{\abs{\Omega}}
   \abs{z^\ast-u^\ast}^2(z^\ast+u^\ast)^{p-2}\,ds\\
   &+\frac pq\int_0^{\abs{\Omega}}
   (u^\ast)^{p-q}\bigl((z^\ast)^q-(u^\ast)^q\bigr)\,ds.
  \end{split}
 \]
 The last integral is non-negative. Indeed, by \eqref{eq:Chiti-majorization}, the function
 \[
  F(s):=\int_0^s\bigl((z^\ast(t))^q-(u^\ast(t))^q\bigr)\,dt
 \]
 is non-negative and vanishes at both endpoints. Since $(u^\ast)^{p-q}$ is non-increasing, integration by parts gives
 \[
  \int_0^{\abs{\Omega}}(u^\ast)^{p-q}F'(s)\,ds
  =-\int_{[0,\abs{\Omega}]}F\,d\bigl((u^\ast)^{p-q}\bigr)
  \ge0.
 \]
which proves \eqref{eq:weighted-Lq}.\medskip

 To obtain \eqref{eq:core-L2-Lq}, we distinguish two cases. If $p\le2$, the uniform bound \eqref{eq:Chiti-uniform-bound} implies
 \[
  \abs{z^\ast-u^\ast}^2(z^\ast+u^\ast)^{p-2}
  \ge(2\norm{z}_\infty)^{p-2}\abs{z^\ast-u^\ast}^2.
 \]
 While, if $p>2$, we use
 \[
  \abs{z^\ast-u^\ast}^p
  \le\abs{z^\ast-u^\ast}^2(z^\ast+u^\ast)^{p-2}
 \]
 and H\"older's inequality to write
 \[
  \int_0^{\abs{B}}\abs{z^\ast-u^\ast}^2\,ds
  \le\abs{B}^{1-2/p}
  \left(\int_0^{\abs{B}}\abs{z^\ast-u^\ast}^p\,ds\right)^{2/p}.
 \]
 In both cases, \eqref{eq:weighted-Lq} gives the desired estimate.\medskip

 Finally, the $L^q$ normalization yields
 \[
  \begin{split}
   \int_{\abs{B}}^{\abs{\Omega}}(u^\ast)^q\,ds
   &=\int_0^{\abs{B}}\bigl((z^\ast)^q-(u^\ast)^q\bigr)\,ds\\
   &\le q\int_0^{\abs{B}}
   \abs{z^\ast-u^\ast}(z^\ast+u^\ast)^{q-1}\,ds.
  \end{split}
 \]
 If $p\le2q$, then
 \[
  \begin{split}
   \abs{z^\ast-u^\ast}(z^\ast+u^\ast)^{q-1}
   &=\left[\abs{z^\ast-u^\ast}^2(z^\ast+u^\ast)^{p-2}\right]^{1/2}
   (z^\ast+u^\ast)^{q-p/2}\\
   &\le C\left[\abs{z^\ast-u^\ast}^2(z^\ast+u^\ast)^{p-2}\right]^{1/2}.
  \end{split}
 \]
 Cauchy--Schwarz and \eqref{eq:weighted-Lq} therefore give a bound by $C\delta^{1/2}$. If $p>2q$, we instead use
 \[
  \begin{split}
   \abs{z^\ast-u^\ast}(z^\ast+u^\ast)^{q-1}
   &=\left[\abs{z^\ast-u^\ast}^2(z^\ast+u^\ast)^{p-2}\right]^{q/p}
   \left(\frac{\abs{z^\ast-u^\ast}}{z^\ast+u^\ast}\right)^{1-2q/p}\\
   &\le\left[\abs{z^\ast-u^\ast}^2(z^\ast+u^\ast)^{p-2}\right]^{q/p}.
  \end{split}
 \]
 H\"older's inequality on $(0,\abs{B})$ and \eqref{eq:weighted-Lq} give a bound by $C\delta^{q/p}$, completing the proof.
\end{proof}

\subsection{Quantitative estimates for the reference set}

As described at the beginning of the section, we set
\begin{equation}\notag
 \bar t:=u^\ast(\abs{B}),\qquad
 U:=\{u>\bar t\},\qquad v:=(u-\bar t)_+.
\end{equation}
The Faber--Krahn inequality gives $\abs{\Omega}\ge\abs{B}$. If equality holds, its rigidity statement already gives the desired conclusions. We may therefore assume that $\abs{\Omega}>\abs{B}$. Since $u$ is a nonconstant real-analytic function on the connected set $\Omega$, every positive level set has measure zero. Consequently, $\bar t>0$ and
\begin{equation}\notag
 \abs{U}=\abs{B},\qquad v\in H_0^1(U).
\end{equation}
We first control the level $\bar t$ and the first eigenvalue of $U$.

\begin{lemma}\label{lem:Chiti-core-eigenvalue}
 There exists $C=C(n,p,q)>0$ such that
 \begin{equation}\label{eq:level-Lq}
  \bar t\le C\delta^{1/(p+1)}.
 \end{equation}
 Moreover, there exists $\delta_0=\delta_0(n,p,q)>0$ such that, if $\delta\le\delta_0$, then
 \begin{equation}\label{eq:core-eigenvalue-Lq}
  0\le\lambda_1(U)-1\le C\bar t
  \le C\delta^{1/(p+1)}.
 \end{equation}
\end{lemma}

\begin{proof}
 The radial eigenfunction $z$ is smooth up to the boundary of the fixed ball $B$ and vanishes there. Hence there exists $K=K(n,q)>0$ such that
 \[
  z^\ast(s)\le K(\abs{B}-s),\qquad 0\le s\le\abs{B}.
 \]
 Also, $\abs{B}\bar t^q\le\norm{u}_q^q=1$. Increasing $K$ if necessary, we may therefore assume that $\bar t/(2K)\le\abs{B}$. For
 \[
  s\in\left(\abs{B}-\frac{\bar t}{2K},\abs{B}\right),
 \]
 we have $u^\ast(s)\ge\bar t$ and $z^\ast(s)\le\bar t/2\le u^\ast(s)/2$. Thus
 \[
  \begin{split}
   \abs{z^\ast-u^\ast}^2(z^\ast+u^\ast)^{p-2}
   &=\left(\frac{u^\ast-z^\ast}{u^\ast+z^\ast}\right)^2
   (u^\ast+z^\ast)^p\\
   &\ge\frac19\bar t^p.
  \end{split}
 \]
 Integrating over this interval and applying \eqref{eq:weighted-Lq}, we obtain
 \[
  \frac{\bar t^{p+1}}{18K}\le C\delta,
 \]
 which proves \eqref{eq:level-Lq}.\medskip

 Before estimating the Rayleigh quotient, we establish a uniform lower bound for $\norm{v}_2$. Since $v^\ast=(u^\ast-\bar t)_+$ and $\abs{U}=\abs{B}$, the reverse triangle inequality gives
 \begin{equation}\label{eq:truncation-norm-Lq}
  \abs{\norm{v}_2-\norm{z}_2}
  \le\norm{u^\ast-z^\ast}_{L^2(0,\abs{B})}
  +\abs{B}^{1/2}\bar t.
 \end{equation}
 By \eqref{eq:core-L2-Lq} and \eqref{eq:level-Lq}, the right-hand side tends to zero with $\delta$. Since $\norm{z}_2>0$ depends only on $n,q$, we can choose $\delta_0=\delta_0(n,p,q)>0$ such that
 \begin{equation}\label{eq:truncation-bounds-Lq}
  \frac{\norm{z}_2}{2}\le\norm{v}_2\le2\norm{z}_2
 \end{equation}
 whenever $\delta\le\delta_0$.\medskip

 Testing $-\Delta u=u$ with $v$ now yields the exact identity
 \begin{equation}\label{eq:Chiti-truncated-energy}
  \int_U\abs{\nabla v}^2\,dx
  =\int_U uv\,dx
  =\norm{v}_2^2+\bar t\int_U v\,dx.
 \end{equation}
 On the other hand, $\abs{U}=\abs{B}$ and the Faber--Krahn inequality imply $\lambda_1(U)\ge1$. Using $v$ in the Rayleigh quotient and then \eqref{eq:truncation-bounds-Lq}, we conclude that
 \[
  0\le\lambda_1(U)-1
  \le\bar t\frac{\displaystyle\int_U v\,dx}{\norm{v}_2^2}
  \le\frac{\abs{B}^{1/2}}{\norm{v}_2}\bar t
  \le C\bar t.
 \]
 Together with \eqref{eq:level-Lq}, this proves \eqref{eq:core-eigenvalue-Lq}.
\end{proof}

We next compare $v$ with the first eigenfunction of the ball. This requires both the spectral gap of $U$ and the quantitative Faber--Krahn estimate for eigenfunctions.

\begin{lemma}\label{prop:truncated-Lq}
 There exist $\delta_0=\delta_0(n,p,q)>0$ and $C=C(n,p,q)>0$ such that, if $\delta\le\delta_0$, then there is $x_\Omega\in\R^n$ satisfying
 \[
  \norm{v-z(\,\cdot-x_\Omega)}_{2}^2
  \le C\delta^{1/(p+1)}.
 \]
\end{lemma}

\begin{proof}
 The Hong--Krahn--Szeg\H{o} inequality \eqref{eq:HKS} and $\abs{U}=\abs{B}$ give $\lambda_2(U)\ge2^{2/n}$. By \eqref{eq:core-eigenvalue-Lq}, decreasing $\delta_0$ if necessary, we obtain
 \begin{equation}\label{eq:gap-Lq}
  \lambda_2(U)-\lambda_1(U)
  \ge\frac12(2^{2/n}-1)>0.
 \end{equation}
 In particular, the first eigenvalue of $U$ is simple.\medskip

 Let $w$ be its non-negative $L^2$-normalized first eigenfunction, and write
 \[
  R(v):=\dfrac{\displaystyle\int_U\abs{\nabla v}^2\,dx}{\displaystyle\int_U v^2\,dx}.
 \]
 Expanding $v/\norm{v}_2$ in an orthonormal basis of Dirichlet eigenfunctions on $U$ gives
 \[
  R(v)-\lambda_1(U)
  \ge\bigl(\lambda_2(U)-\lambda_1(U)\bigr)
  \left[1-\left(\int_U\frac{v}{\norm{v}_2}w\,dx\right)^2\right].
 \]
 Since $v,w\ge0$ and both normalized functions have unit $L^2$ norm,
 \[
  \begin{split}
   \norm{\frac{v}{\norm{v}_2}-w}_2^2
   &=2\left(1-\int_U\frac{v}{\norm{v}_2}w\,dx\right)\\
   &\le2\left[1-\left(\int_U\frac{v}{\norm{v}_2}w\,dx\right)^2\right].
  \end{split}
 \]
 Moreover, \eqref{eq:Chiti-truncated-energy} and \eqref{eq:truncation-bounds-Lq} imply
 \[
  0\le R(v)-\lambda_1(U)\le R(v)-1\le C\bar t.
 \]
 Combining these estimates with \eqref{eq:gap-Lq}, we find
 \[
  \norm{\frac{v}{\norm{v}_2}-w}_2^2\le C\bar t.
 \]

\noindent By \autoref{thm:eigenfunction:fk}, there exists $x_\Omega\in\R^n$ such that
 \[
  \norm{w-\frac{z(\,\cdot-x_\Omega)}{\norm{z}_2}}_2^2
  \le C\bigl(\lambda_1(U)-1\bigr)
  \le C\bar t.
 \]
 Therefore
 \begin{equation}\label{eq:v-z-normalized-Lq}
  \norm{\frac{v}{\norm{v}_2}
  -\frac{z(\,\cdot-x_\Omega)}{\norm{z}_2}}_2^2
  \le C\bar t.
 \end{equation}
 Finally, \eqref{eq:truncation-norm-Lq}, \eqref{eq:truncation-bounds-Lq}, and \eqref{eq:v-z-normalized-Lq} yield
 \[
  \begin{split}
   \norm{v-z(\,\cdot-x_\Omega)}_2^2
   &\le2\norm{v}_2^2
   \norm{\frac{v}{\norm{v}_2}
   -\frac{z(\,\cdot-x_\Omega)}{\norm{z}_2}}_2^2
   +2\abs{\norm{v}_2-\norm{z}_2}^2\\
   &\le C\left(\bar t+\delta^{\min\{1,2/p\}}+\bar t^2\right)\\
   &\le C\delta^{1/(p+1)},
  \end{split}
 \]
 where we used \eqref{eq:level-Lq}, $\delta_0\le1$, and $\min\{1,2/p\}>1/(p+1)$.
\end{proof}

\subsection{Proof of \autoref{thm:intro-Chiti}}

It remains to estimate the truncation error and to pass from $L^2$ to $L^q$. We record the former estimate for every finite $r\ge q$, since it will also be used in the final corollary.

\begin{lemma}
 For every finite $r\ge q$, there exists $C=C(n,p,q,r)>0$ such that
 \begin{equation}\label{eq:discarded-r-Lq}
  \norm{u-v}_r^r
  \le\abs{B}\bar t^r
  +C\bar t^{r-q}\delta^{\min\{1/2,q/p\}}.
 \end{equation}
 In particular, if $\delta\le1$, then
 \begin{equation}\label{eq:Chiti-truncation-error-power}
  \norm{u-v}_r^{\max\{2,r\}}
  \le C\delta^{1/(p+1)}.
 \end{equation}
\end{lemma}

\begin{proof}
 Since $u-v=\min\{u,\bar t\}$, equimeasurability gives
 \[
  \begin{split}
   \norm{u-v}_r^r
   &=\abs{B}\bar t^r
   +\int_{\abs{B}}^{\abs{\Omega}}(u^\ast)^r\,ds\\
   &\le\abs{B}\bar t^r
   +\bar t^{r-q}\int_{\abs{B}}^{\abs{\Omega}}(u^\ast)^q\,ds.
  \end{split}
 \]
 \eqref{lemma3.3last} proves \eqref{eq:discarded-r-Lq}.\medskip

 If $\delta\le1$, then \eqref{eq:level-Lq} bounds $\bar t$ uniformly. Consequently,
 \[
  \norm{u-v}_r^{\max\{2,r\}}
  \le C\left(
  \bar t^{\max\{2,r\}}
  +\delta^{\min\{1/2,q/p\}\max\{2,r\}/r}
  \right).
 \]
 Since
 \[
  \min\{1/2,q/p\}>\frac1{p+1},\qquad
  \frac{\max\{2,r\}}r\ge1,
 \]
 another application of \eqref{eq:level-Lq} gives \eqref{eq:Chiti-truncation-error-power}.
\end{proof}

We can now conclude the proof of the main theorem.

\begin{proof}[Proof of \autoref{thm:intro-Chiti}]
 First assume \eqref{eq:Chiti-normalization} and $\delta\le\delta_0\le1$, where $\delta_0$ is chosen as in \autoref{prop:truncated-Lq}. Let $x_\Omega$ be the point provided by that lemma. The function $v-z(\,\cdot-x_\Omega)$ is supported on a set of measure at most $2\abs{B}$ and, by \eqref{eq:Chiti-uniform-bound}, its $L^\infty$ norm is at most $2\norm{z}_\infty$. Thus \autoref{lem:change-norm-Lq} gives
 \[
  \norm{v-z(\,\cdot-x_\Omega)}_q^{\max\{2,q\}}
  \le C\norm{v-z(\,\cdot-x_\Omega)}_2^2
  \le C\delta^{1/(p+1)}.
 \]
 Together with \eqref{eq:Chiti-truncation-error-power} for $r=q$ and the triangle inequality, this yields
 \begin{equation}\label{eq:Chiti-quantitative-normalized}
  \inf_{x_0\in\R^n}
  \norm{u-z(\,\cdot-x_0)}_q^{\max\{2,q\}}
  \le C\delta^{1/(p+1)}.
 \end{equation}
 If $\delta>\delta_0$, then the $L^q$ normalization gives $\norm{u-z(\,\cdot-x_0)}_q\le2$ for every $x_0$. Increasing $C$ if necessary, \eqref{eq:Chiti-quantitative-normalized} therefore holds for every deficit.\medskip

 Finally, let $\Omega$ have arbitrary first eigenvalue $\lambda>0$, and let $u,z$ satisfy \eqref{eq:q-normalization}. Set
 \begin{equation}\notag
  \widehat\Omega:=\sqrt{\lambda}\,\Omega,\qquad
  \widehat u(x):=\lambda^{-\frac{n}{2q}}u(x/\sqrt{\lambda}),\qquad
  \widehat z(x):=\lambda^{-\frac{n}{2q}}z(x/\sqrt{\lambda}).
 \end{equation}
 Then $\lambda_1(\widehat\Omega)=1$, $\sqrt{\lambda}\,B^\lambda=B$, and $\widehat u,\widehat z$ have unit $L^q$ norm. A change of variables gives
 \[
  \begin{split}
   \norm{\widehat z}_{L^p(B)}^p
   -\norm{\widehat u}_{L^p(\widehat\Omega)}^p
   &=\lambda^{-\frac{n(p-q)}{2q}}
   \left(\norm{z}_{L^p(B^\lambda)}^p
   -\norm{u}_{L^p(\Omega)}^p\right)\\
   &=\delta_{\mathrm{Ch}}(\Omega).
  \end{split}
 \]
 Moreover, for every $x_0\in\R^n$,
 \[
  \norm{\widehat u-\widehat z(\,\cdot-\sqrt{\lambda}\,x_0)}_q
  =\norm{u-z(\,\cdot-x_0)}_q.
 \]
 Applying \eqref{eq:Chiti-quantitative-normalized} on $\widehat\Omega$ proves the scale-invariant estimate
 \begin{equation}\label{eq:Chiti-quantitative-invariant}
  \delta_{\mathrm{Ch}}(\Omega)
  \ge c(n,p,q)
  \left(\inf_{x_0\in\R^n}
  \norm{u-z(\,\cdot-x_0)}_q^{\max\{2,q\}}\right)^{p+1},
 \end{equation}
 which is the assertion using that $\dfrac{u}{\norm{u}_{q}}$ and $\dfrac{z}{\norm{z}_{q}}$satisfies the normalization \eqref{eq:q-normalization}.
\end{proof}

We now consider the disconnected case. Since in this case the first eigenfunction could not be unique up to multiplication by a constant, the Chiti deficit may depend on the chosen eigenfunction. Hence, for every non-negative first Dirichlet eigenfunction $u$ on $\Omega$, we define
\[
\delta_{\mathrm{Ch}}(\Omega,u)
:=
\lambda_1(\Omega)^{\frac{n(q-p)}{2q}}
\left(
\frac{\|z\|_{L^p(B^\lambda)}^p}{\|z\|_{L^q(B^\lambda)}^p}
-
\frac{\|u\|_{L^p(\Omega)}^p}{\|u\|_{L^q(\Omega)}^p}
\right),
\]
where $B^\lambda$ is the ball satisfying
$\lambda_1(B^\lambda)=\lambda_1(\Omega)$ and $z$ is any positive first Dirichlet eigenfunction on $B^\lambda$.

\begin{corollary}\label{cor:Chiti-disconnected}
Let $1\le q<p<\infty$. There exists a positive constant
$C=C(n,p,q)$ such that, for every bounded open set
$\Omega\subset\mathbb R^n$ and every non-negative first Dirichlet
eigenfunction $u$ on $\Omega$,
\[
\delta_{\mathrm{Ch}}(\Omega,u)
\ge C
\left(
\inf_{x_0\in\mathbb R^n}
\left\|
\frac{u}{\|u\|_q}
-
\frac{z(\,\cdot-x_0)}{\|z\|_q}
\right\|_q^{\max\{2,q\}}
\right)^{p+1},
\]
where $B^\lambda$ is the ball satisfying
$\lambda_1(B^\lambda)=\lambda_1(\Omega)$ and $z$ is any positive
first Dirichlet eigenfunction on $B^\lambda$.
\end{corollary}

\begin{proof} Let $u$ be any non-negative first eigenfunction.
By scaling and normalization, we may assume, without loss of generality, that
\[
    \lambda_1(\Omega)=1,\qquad\norm{u}_{L^q(\Omega)}=1.
\]
We also normalize the first eigenfunction $z$ on $B^1$ so that $\norm{z}_{L^q(B^1)}=1$.\medskip

Let $\{\Omega_i\}_{i\in I}$ be the connected components of $\Omega$ on which $u$ is not identically zero. Since $u$ is a non-negative first eigenfunction, it satisfies
\[
    \lambda_1(\Omega_i)=\lambda_1(\Omega)=1.
\]
For every $i\in I$, define
\[
    t_i:=\int_{\Omega_i}u^q\,dx,\qquad u_i:=\frac{u\,\mathbf1_{\Omega_i}}{t_i^{1/q}}.
\]
Clearly, $t_i>0$, $\sum_{i\in I}t_i=1$, and $\norm{u_i}_{q}=1$.\medskip

By the definition of the Chiti deficit and the disjointness of the components, we have
\begin{equation}\label{eq:deficit_disjoint}
    \begin{aligned}
        \delta_{\mathrm{Ch}}(\Omega,u)&=\norm{z}_{p}^p-\norm{u}_{p}^p\\
        &=\norm{z}_{p}^p-\sum_{i\in I}t_i^{p/q}\norm{u_i}_{p}^p\\
        &=\norm{z}_{p}^p\left(1-\sum_{i\in I}t_i^{p/q}\right)+\sum_{i\in I}t_i^{p/q}\delta_{\mathrm{Ch}}(\Omega_i).
     \end{aligned}
\end{equation}
Where we used the fact that all components have the same first eigenvalue. Without loss of generality we can assume that $\max_{i\in I}t_i=t_1$, thus we have
\[
    \sum_{i\in I}t_i^{p/q}\leq t_1^{p/q-1}.
\]
Therefore, using the non-negativity of $\delta_{\mathrm{Ch}}(\Omega_i)$ in \eqref{eq:deficit_disjoint},
\begin{equation}\label{eq:chiti_mass_gap}
    \delta_{\mathrm{Ch}}(\Omega,u)\geq\norm{z}_{p}^p\left(1-\sum_{i\in I}t_i^{p/q}\right)\geq\norm{z}_{p}^p\left(1-t_1^{p/q-1}\right)\geq c(n,p,q)(1-t_1).
\end{equation}
Moreover, identity \eqref{eq:deficit_disjoint} immediately yields
\begin{equation}\label{eq:deficit_component}
    t_1^{p/q}\delta_{\mathrm{Ch}}(\Omega_1)\leq\delta_{\mathrm{Ch}}(\Omega,u).
\end{equation}
We now apply \autoref{thm:intro-Chiti} to the connected component $\Omega_1$, that is
\begin{equation}\label{eq:chiti_on_omega1}
    \inf_{x_0\in\R^n}\norm{u_1-z(\,\cdot-x_0)}_{q}^{\max\{2,q\}}\leq C(n,p,q)\delta_{\mathrm{Ch}}(\Omega_1)^{1/(p+1)}.
\end{equation}
On the other hand, for every $x_0\in\R^n$, the triangle inequality gives
\begin{equation*}
    \begin{aligned}
        \norm{u-z(\,\cdot-x_0)}_{q}&\leq\norm{u-t_1^{1/q}u_1}_{q}+t_1^{1/q}\norm{u_1-z(\,\cdot-x_0)}_{q}+(1-t_1^{1/q})\norm{z}_{q}\\
        &=(1-t_1)^{1/q}+t_1^{1/q}\norm{u_1-z(\,\cdot-x_0)}_{q} +(1-t_1^{1/q})\\
        &\leq t_1^{1/q}\norm{u_1-z(\,\cdot-x_0)}_{q}+2(1-t_1)^{1/q}.
    \end{aligned}
\end{equation*}
Taking the infimum over $x_0$ and using \eqref{eq:chiti_on_omega1}, we deduce
\begin{equation*}
    \inf_{x_0\in\R^n}\norm{u-z(\,\cdot-x_0)}_{q}^{\max\{2,q\}}\leq C(n,p,q)\left(t_1^{\max\{2,q\}/q}\delta_{\mathrm{Ch}}(\Omega_1)^{1/(p+1)}+(1-t_1)^{\max\{2,q\}/q}\right).
\end{equation*}
By \eqref{eq:deficit_component}, since $\max\{2,q\}>p/(p+1)$ and $0<t_1\leq1$, we have
\[
    t_1^{\max\{2,q\}/q}\delta_{\mathrm{Ch}}(\Omega_1)^{1/(p+1)}\leq t_1^{\frac1q(\max\{2,q\}-\frac{p}{p+1})}\delta_{\mathrm{Ch}}(\Omega,u)^{1/(p+1)}\leq\delta_{\mathrm{Ch}}(\Omega,u)^{1/(p+1)}.
\]
Furthermore, since $\max\{2,q\}/q\geq1$, estimate \eqref{eq:chiti_mass_gap} implies
\[
    (1-t_1)^{\max\{2,q\}/q}\leq (1-t_1)^{1/(p+1)}\leq C(n,p,q)\delta_{\mathrm{Ch}}(\Omega,u)^{1/(p+1)}.
\]
Combining the previous estimates, we obtain
\[
    \inf_{x_0\in\R^n}\norm{u-z(\,\cdot-x_0)}_{L^q(\R^n)}^{\max\{2,q\}}\leq C(n,p,q)\delta_{\mathrm{Ch}}(\Omega,u)^{1/(p+1)},
\]
hence the assertion.
\end{proof}

\subsection{Further consequences}

We conclude with a few consequences of the preceding estimates.

\begin{corollary}[Almost radial symmetry]
 Under the assumptions of \autoref{thm:intro-Chiti}, let $u^\sharp$ denote the Schwarz rearrangement of $u$. Then
 \[
  \delta_{\mathrm{Ch}}(\Omega)\geq C\left(\inf_{x_0\in\R^n}
  \norm{\frac{u-u^\sharp(\,\cdot-x_0)}{\norm{u}_{q}}}_{L^q(\R^n)}^{\max\{2,q\}}\right)^{(p+1)}.
 \]
\end{corollary}

\begin{proof}
 Assuming for simplicity that $\norm{u}_{q}=1$, since Schwarz rearrangement is a contraction in $L^q$ and the rearrangement of $z(\,\cdot-x_0)$ is $z$, for every $x_0\in\R^n$ we have
 \[
  \norm{u^\sharp-z}_q
  \le\norm{u-z(\,\cdot-x_0)}_q.
 \]
 Therefore
 \[
  \begin{split}
   \norm{u-u^\sharp(\,\cdot-x_0)}_q
   &\le\norm{u-z(\,\cdot-x_0)}_q+\norm{z-u^\sharp}_q\\
   &\le2\norm{u-z(\,\cdot-x_0)}_q.
  \end{split}
 \]
 Taking the infimum and applying \eqref{eq:Chiti-quantitative-invariant} proves the assertion.
\end{proof}

\begin{corollary}
Under the assumptions of \autoref{thm:intro-Chiti},
for every finite $r\ge q$ there exists a positive constant
$C=C(n,p,q,r)$ such that
\begin{equation}\notag
    \delta_{\mathrm{Ch}}(\Omega)\ge C\left(\inf_{x_0\in\R^n}\left[\lambda_1(\Omega)^{\frac n2(\frac1r-\frac1q)}\norm{\frac{u}{\norm{u}_{q}}-\frac{z(\,\cdot-x_0)}{\norm{z}_{q}}}_{r}\right]^{\max\{2,r\}}\right)^{p+1}.
\end{equation}
\end{corollary}
\begin{proof}
We first prove the assertion under the normalizations
\[
    \lambda=1,\qquad\norm{u}_{q}=\norm{z}_{q}=1.
\]
Assume first that $\delta\le\delta_0\le1$, where $\delta_0>0$ is the smallness threshold appearing in the proof of \autoref{thm:intro-Chiti}.

Let $v$ be the truncated function constructed in that proof, and let $x_\Omega$ be the corresponding translation parameter. The uniform support and $L^\infty$ bounds, together with \autoref{lem:change-norm-Lq} and \autoref{prop:truncated-Lq}, yield
\[
    \norm{v-z(\,\cdot-x_\Omega)}_{2}^2\le C\delta^{1/(p+1)}.
\]
Moreover, the same uniform bounds imply, for every finite $r\ge q$,
\[
    \norm{v-z(\,\cdot-x_\Omega)}_{r}^{\max\{2,r\}}\le C\norm{v-z(\,\cdot-x_\Omega)}_{2}^2.
\]
Indeed, for $r\le2$ this follows from the uniform bound on the measure of the supports, whereas for $r\ge2$ it follows from the uniform $L^\infty$ bound. Consequently,
\begin{equation}\label{eq:Chiti-higher-norms-truncated}
    \norm{v-z(\,\cdot-x_\Omega)}_{r}^{\max\{2,r\}}\le C\delta^{1/(p+1)}.
\end{equation}

On the other hand, the truncation error estimate \eqref{eq:Chiti-truncation-error-power}, combined with the uniform $L^\infty$ bound and $r\ge q$, gives
\begin{equation}\label{eq:Chiti-higher-norms-error}
    \norm{u-v}_{r}^{\max\{2,r\}}\le C\delta^{1/(p+1)}.
\end{equation}
Therefore, the triangle inequality and \eqref{eq:Chiti-higher-norms-truncated}--\eqref{eq:Chiti-higher-norms-error} imply
\[
    \norm{u-z(\,\cdot-x_\Omega)}_{r}^{\max\{2,r\}}\le C\left(\norm{u-v}_{r}^{\max\{2,r\}}+\norm{v-z(\,\cdot-x_\Omega)}_{r}^{\max\{2,r\}}\right)\le C\delta^{1/(p+1)}.
\]
This proves the statement for sufficiently small deficits.

If $\delta>\delta_0$, the uniform $L^\infty$ estimates and the $L^q$ normalizations give
\[
    \norm{u}_{r}^r\le\norm{u}_{L^\infty}^{r-q}\norm{u}_{q}^q\le C,
\]
and the same estimate holds for $z$. Hence, for every $x_0\in\R^n$,
\[
    \norm{u-z(\,\cdot-x_0)}_{r}^{\max\{2,r\}}\le C\le C\delta_0^{-1/(p+1)}\delta^{1/(p+1)}.
\]
After increasing the constant, we obtain the statement for every $\delta\ge0$.\medskip

Finally, we compute the scaling of the $L^r$ distance. For any $x_0\in\R^n$, setting $\widehat x_0=\sqrt{\lambda}\,x_0$ and $\widehat u$ for the normalized eigenfunction, we have
\begin{equation*}
    \norm{\widehat u-\widehat z(\,\cdot-\widehat x_0)}_{r}^r=\lambda^{-\frac{nr}{2q}}\int_{\R^n}\abs{\frac{u(x/\sqrt{\lambda})}{\norm{u}_{q}}-\frac{z(x/\sqrt{\lambda}-x_0)}{\norm{z}_{q}}}^r\,dx=\lambda^{\frac n2-\frac{nr}{2q}}\int_{\R^n}
    \abs{\frac{u(y)}{\norm{u}_{q}}-
    \frac{z(y-x_0)}{\norm{z}_{q}}}^r\,dy.
\end{equation*}
Applying the previously obtained result to $\widehat u$ and $\widehat z$, and using the invariance of the deficit together with the $L^r$ scaling computed above, we obtain
\[
    C\inf_{x_0\in\R^n}\left[\lambda^{\frac n2(\frac1r-\frac1q)}\norm{\frac{u}{\norm{u}_{q}}-\frac{z(\,\cdot-x_0)}{\norm{z}_{q}}}_{r}\right]^{\max\{2,r\}}\le \delta_{\mathrm{Ch}}(\Omega)^{1/(p+1)}.\qedhere
\]
\end{proof}

\section{Quantitative Kohler--Jobin inequality }

In this section, we prove \autoref{thm:intro-KJ}. Let $\Omega\subset\R^n$ be a bounded, open, and connected set and let $u>0$ be the first eigenfunction on $\Omega$ normalized by
\[
 \norm{u}_{L^2(\Omega)}=1.
\]
For simplicity, throughout the proof we assume that
\begin{equation}\notag
 T(\Omega)=1,
\end{equation}
so that, letting $r_0>0$ be the radius of the ball with unitary torsional rigidity, the Kohler--Jobin deficit is exactly
\begin{equation}\notag
\begin{split}\delta=\delta_{\mathrm{KJ}}(\Omega)&=T(\Omega)^\frac{2}{n+2}\lambda_1(\Omega)-T(B_{r_0})^\frac{2}{n+2}\lambda_1(B_{r_0})\\[5pt]&=\lambda_1(\Omega)-\lambda_1(B_{r_0}).\end{split}
\end{equation}\medskip

The section is structured as follows. We will start by recalling the definition of modified torsional rigidity and of the so-called Kohler--Jobin rearrangement of $u$, which is defined on a ball $B_r$, of smaller volume, such that
\[T(B_r)=\Tmod(\Omega;u),\]
and show that 
\[0\le r_0-r\lesssim \delta.\]
Then, as in the proof of \autoref{thm:intro-Chiti}, we will introduce the reference set $U=\set{u>\bar{t}}$ such that $\abs{U}=\abs{B_r}$ and prove that, if  $\delta$ is sufficiently small, then 
\[\bar{t}\lesssim \delta^{1/2}.\]
From the previous inequality we will eventually deduce that
\[\lambda_1(U)-\lambda_1(B_r)\lesssim \delta^{1/2},\]
so that, applying \eqref{eq:AKN-res}, we obtain
\[\norm{\mathcal{R}_U-\mathcal{R}_{B_r(x_\Omega)}}_{L^\infty\to L^2}\lesssim\delta^{\frac{1}{4}},\]
for some $x_\Omega\in R^n$. Hence, we will conclude the proof by estimating the remainder terms
\[\norm{\mathcal{R}_\Omega-\mathcal{R}_{U}}_{L^\infty\to L^2},\qquad\text{and}\qquad\lVert\mathcal{R}_{B_{r}}-\mathcal{R}_{B_{r_0}}\rVert_{L^\infty\to L^2}.\]\medskip

In the notation of the sketch of the proof of \autoref{KJ-rearrangment}, let $\Phi\,\colon [0,M]\to[0,+\infty)$ be such that the torsional rigidity of the ball $\set{v>\Phi(t)}$ is exactly $\mathcal{T}(t)$, and denote by $\beta(t)$ its measure. Then, $\Phi(t)\ge t$, and letting
\[D(t):=\Phi(t)-t,\]
we have 
\[D'(t)=\dfrac{\mu(t)}{\beta(t)}-1,\]
and  
\[\int_{B_r} v^2\,dx-1=\int_{B_r} v^2\,dx-\int_\Omega u^2\,dx=2\int_0^M D(t)\mu(t)\,dt.\]
We have the following lemma.

\begin{lemma}\label{lem:KJ-exact}
There exists a constant $C_n=C(n)>0$ such that

\begin{equation}\label{eq:KJ-exact}
 2\int_0^M D(t)\mu(t)\,dt=\int_{B_r}v^2\,dx-1\le C_n \delta,
\end{equation}
and
\begin{equation}\label{eq:r-close}
 0\le r_0-r\le C_n \delta.
\end{equation}
\end{lemma}

\begin{proof}
By definition, $B_r$ is the ball such that 
\[T(B_r)=\Tmod(\Omega;u)\le T(\Omega)=1=T(B_{r_0}),\]
so that we have $r\le r_0$.
Moreover, recall that, with our normalization $T(\Omega)=1$, the Kohler--Jobin deficit is exactly 
\[\delta = \lambda_1(\Omega)-\lambda_1(B_{r_0}).\]
Hence, by the properties of the rearrangement,
\[r^{-2}\lambda_1(B_1)=\lambda_1(B_r)\le \dfrac{\lambda_1(\Omega)}{\int_{B_r} v^2\,dx}=\dfrac{r_0^{-2}\,\lambda_1(B_1)+\delta}{\int_{B_r} v^2\,dx},\]
that is
\[\int_{B_r}v^2\,dx\le \left(\dfrac{r}{r_0}\right)^{2}+r^2\dfrac{\delta}{\lambda_1(B_1)}.\]
Using the fact that $r\le r_0$, we obtain
\[\int_{B_r}v^2\,dx-1\le \dfrac{\delta}{\lambda_1(B_{r_0})},\]
hence \eqref{eq:KJ-exact}. On the other hand, using the inequality 
\[\int_{B_r} v^2\,dx\ge 1,\]
we have
\[r_0-r\le\dfrac{r^2 \delta}{\lambda_1(B_{r_0})(r_0+r)}\le\dfrac{r_0 \delta}{2 \lambda_1(B_{r_0})},\]
that is \eqref{eq:r-close}.
\end{proof}

\subsection{Quantitative estimates for the reference set}
As in the proof of \autoref{thm:intro-Chiti}, we now introduce the reference set and prove that its distance, in the resolvent sense, from the ball $B_r$ can be estimated in terms of the Kohler--Jobin deficit of $\Omega$.

Set
\begin{equation}\notag
 \bar t:=u^\ast(\abs{B_r}),\qquad U:=\{u>\bar t\}.
\end{equation}
If $\abs{B_r}=|\Omega|$, from 
\[T(B_r)=\Tmod(\Omega;u)\le T(\Omega)\]
and by the rigidity of the Saint--Venant inequality, we have that $\Omega$ is a ball.  Otherwise $\bar t>0$ and the analyticity of $u$ gives
\begin{equation}\notag
 |U|=|B_r|.
\end{equation}
We prove the following proposition
\begin{proposition}\label{prop:torsion-core-gap}
    There exist constants $C_n=C(n)>0$ and $\delta_0=\delta_0(n)>0$,  and a point $x_\Omega\in\R^n$, such that, if $\delta<\delta_0$, then
    \begin{equation}\notag
        \norm{\mathcal{R}_U-\mathcal{R}_{B_r(x_\Omega)}}_{L^\infty\to L^2}\le C_n \delta^{\frac{1}{4}}.
    \end{equation}
\end{proposition}

\begin{remark}\label{rem:boundnorm}
We notice that, if $\delta$ is sufficiently small, both the $L^1$ and the $L^\infty$ norms of $u$ are uniformly bounded in $n$. Indeed, using $u$ as a test function in the definition of the torsional rigidity, we have 
\begin{equation}\notag
 \left(\int_\Omega u\,dx\right)^2\le T(\Omega)\int_\Omega\abs{\nabla u}^2
 =\lambda_1(\Omega)=\lambda_1(B_{r_0})+\delta.
\end{equation}
Similarly, by Chiti's reverse H\"odler inequality, we have that 
\begin{equation}\notag
\norm u_\infty\le C_n\lambda_1(\Omega)^{n/4}\norm u_2=C_n \left(\lambda_1(B_{r_0})+\delta\right)^{n/4}.\end{equation}
Hence, if $\delta\le1$, we have $\norm{u}_1,\norm{u}_\infty\le C_n.$
\end{remark}

We now prove that, if $\delta$ is sufficiently small, the level $\bar{t}$ can be controlled with the Kohler--Jobin deficit.
\begin{proposition}\label{prop:tbar-KJ}
There exist $\delta_0(n)>0$ and $C_n>0$ such that, if $\delta\le\delta_0$, then
\begin{equation}\notag
 \bar t\le C_n\delta^{1/2}.
\end{equation}
\end{proposition}

\begin{proof}
The idea of the proof is to show that, if $\delta$ is sufficiently small, then 
\[2\int_{\bar{t}}^M D(s)\mu(s)\,ds\ge C_n \bar{t}^2,\]
so that, by the estimate \eqref{eq:KJ-exact}, we get the assertion.\medskip

We start by proving that there exist a dimensional constant $C_n$, such that, if $t\le t_0=t_0(n)$ and $\delta\le1$, then
\begin{equation}\label{massestimate}\int_t^M \mu(s)\,ds\ge C_n.\end{equation}
Indeed 
\[\norm{u}_\infty\int_t^M \mu(s)\,ds=\norm{u}_\infty\int_\Omega (u-t)_+\,dx\ge\int_\Omega (u-t)_+^2\,dx\ge\int_\Omega (u^2-2tu)\,dx.\]
If $\delta\le 1$, by \autoref{rem:boundnorm}, the $L^1$ and $L^\infty$ norms of $u$ are uniformly bounded from above, so that
\[\int_t^M \mu(s)\,ds\ge \dfrac{1-C_n t}{C_n}.\]
Then, if $t<t_0=(2C_n)^{-1}$, we get \eqref{massestimate}.\medskip

On the other hand, there exists a dimensional constant $C_n>0$, such that, if $\delta$ is sufficiently small, we have  
\begin{equation}\label{Dgrowth}D(t)\ge C_n t^2\end{equation}
for every $t\le \bar{t}$. Indeed, for every $0\le t\le \bar{t}$, $\mu(t)\ge \mu(\bar{t})=\abs{B_r}$. By \eqref{eq:r-close} 
\[r\ge r_0-C_n \delta,\]
so that if $\delta$ is sufficiently small, for instance $\delta\le \delta_0=r_0(2 C_n)^{-1}$, then $r$ is uniformly bounded from below, hence so are $\mu(t)\ge\abs{B_r}=\omega_n r^n$ and \[\mathcal{T}(0)=T(B_r)=\left(\dfrac{r}{r_0}\right)^{n+2}.\] Recalling \eqref{eq:Tt-formula}, we have that 
\[\mathcal{T}(t)=\Tmod(\Omega; (u-t)_+)=\int_t^M \dfrac{\mu(s)^2}{\lambda_1(\Omega) F(s)}\,ds,\]
where 
\[F(s)=\int_{\Omega_s} u_s\,dx\le \norm{u}_1\le C_n,\]
if $\delta\le 1$. Hence, we have 
\[\mathcal{T}(0)-\mathcal{T}(t)=\int_0^t \dfrac{\mu(s)^2}{\lambda_1(\Omega) F(s)}\,ds\ge C_n t. \]

Recall that, by construction, 
\[D'(t)=\dfrac{\mu(t)-\beta(t)}{\beta(t)},\]
where \[\beta(t)=\omega_n r_0^n\mathcal{T}(t)^\frac{n}{n+2}\] is the measure of the ball whose torsional rigidity is 
$\mathcal{T}(t)$.
Then, for every $t\le\bar{t}$, we have $\beta(t)\le \beta(0) = \abs{B_r}$, so that
\begin{equation}\label{D'growth}D'(t)\ge\dfrac{\omega_n r_o^n\left(\mathcal{T}(0)^\frac{n}{n+2}-\mathcal{T}(t)^\frac{n}{n+2}\right)}{\abs{B_r}}\ge C_n \left(\mathcal{T}(0)-\mathcal{T}(t)\right)\ge C_n t,\end{equation}
where we also used that $\mathcal{T}(0)\le1$. Integrating \eqref{D'growth} from $0$ to $t$, recalling that $D(0)=0$, we have \eqref{Dgrowth}.\medskip

Finally, let $\delta$ be sufficiently small for \eqref{massestimate} and \eqref{Dgrowth} to hold. Then, if, $\bar t\le t_0$, by \eqref{eq:KJ-exact},  we have
\[
C_n \delta \ge 2 \int_t ^M D(s)\mu(s)\,ds
 \ge2D(\bar t)\int_{\bar t}^M\mu(s)\,ds
 \ge C'_n\bar t^2,
\]
that is the assertion.
 On the other hand, if $\bar t>t_0$ we get a contradiction. Indeed, the same estimate gives 
 \[C_n \delta\ge2D( t_0)\int_{ t_0}^M\mu(s)\,ds
 \ge C_n' t_0^2\] which is a contradiction for $\delta$ sufficiently small.
\end{proof}

The previous proposition immediately yields the following estimate.

\begin{proposition}\label{eigclose}
There exist $\delta_0(n)>0$ and $C_n>0$ such that, if $\delta\le\delta_0$, then
\begin{equation}\label{eq:UOmegaeigcolose}
 \lambda_1(U)-\lambda_1(\Omega)\le C_n\delta^{1/2},
\end{equation}
moreover
\begin{equation}\label{eq:UBreigcolose}
 \lambda_1(U)-\lambda_1(B_r)\le C_n\delta^{1/2},
\end{equation}
\end{proposition}
\begin{proof}
    Using $(u-\bar{t})_+$ in the equation $-\Delta u =\lambda_1(\Omega)u$, we have 
\[\int_U \abs{\nabla (u-\bar{t}\,)}^2\,dx = \lambda_1(\Omega) \left(\int_U (u-\bar t\,)^2\,dx+\bar t\int_U (u-\bar t\,)\,dx\right).\]
If $\delta\le1$, by \autoref{rem:boundnorm}, then
\[\int_U (u-\bar t\,)\,dx\le \norm{u}_1\le C_n,\]
and
\[\int_U (u-\bar t\,)^2\,dx\ge \int_\Omega u^2 - 2\bar t\int_U u\ge 1-2\bar{t}C_n.\]
Estimating the Rayleigh quotient of $U$ with $(u-t)_+$ gives
\[\lambda_1(U)\le \dfrac{\displaystyle \int_U \abs{\nabla (u-\bar{t}\,)}^2\,dx}{\displaystyle\int_U (u-\bar t\,)^2\,dx}\le\lambda_1(\Omega)\left(1+\dfrac{\bar t\, C_n}{1-\bar{t}\,C_n}\right).\]
By \autoref{prop:tbar-KJ}, if $\delta$ is sufficiently small, then $\bar{t}\le C_n \delta^{1/2}$, hence
\[\lambda_1(U)-\lambda_1(\Omega)\le C_n \delta^{1/2}.\]
To prove \eqref{eq:UBreigcolose}, we simply notice that, as $r\le r_0$, then $\lambda_1(B_r)\ge\lambda_1(B_{r_0})$, and
\[\lambda_1(\Omega)-\lambda_1(B_r)\le \lambda_1(\Omega)-\lambda_1(B_{r_0})=\delta.\]
Thus, by \eqref{eq:UOmegaeigcolose}, 
\[\lambda_1(U)-\lambda_1(B_r)\le\lambda_1(\Omega)-\lambda_1(B_r)+C_n\delta^{1/2}\le C_n\delta^{1/2}.\]
\end{proof}

\autoref{prop:torsion-core-gap} is now a direct result of \autoref{eigclose} and \autoref{AKN-res}. 
\begin{proof}[proof of \autoref{prop:torsion-core-gap}]
As $\abs{U}=\abs{B_r}$, we know that, by \eqref{eq:AKN-res}, there exists $x_\Omega\in\R^n$ and a dimensional constant $C_n>0$ such that 
\[\norm{\mathcal{R}_U - \mathcal{R}_{B_r(x_\Omega)} }_{L^\infty\to L^2}^2\le C_n \left(\lambda_1(U)-\lambda_1(B_r)\right).\]
Hence, for $\delta$ sufficiently small, estimate \eqref{eq:UBreigcolose} gives the assertion.
\end{proof}

\subsection{Proof of \autoref{thm:intro-KJ}}

As stated in the introduction of the section, we are now left to estimate the remainder terms
\[\norm{\mathcal{R}_\Omega-\mathcal{R}_{U}}_{L^\infty\to L^2},\qquad\text{and}\qquad\lVert\mathcal{R}_{B_{r}}-\mathcal{R}_{B_{r_0}}\rVert_{L^\infty\to L^2}.\]
To do that, we prove the following simple lemma.

\begin{lemma}\label{lem:nested-resolvents}
Let $D\subset E$ be bounded open sets.  Then
\begin{equation}\label{eq:nested-resolvents}
 \norm{\res_E-\res_D}_{L^\infty\to L^2}^2
 \le\frac{T(E)-T(D)}{\lambda_1(E)}.
\end{equation}
\end{lemma}

\begin{proof}
Let $w_E$ and $w_D$ be the torsion functions, extended by zero outside the corresponding domains. By the maximum principle, the operator $\res_E-\res_D$ is positive.  Hence, for every $\norm f_\infty\le1$,
\[
 |(\res_E-\res_D)f|
 \le(\res_E-\res_D)|f|
 \le (\res_E-\res_D)1 = w_E-w_D.
\]
so that
\begin{equation}\label{restorsion}
 \norm{\res_E-\res_D}_{L^\infty\to L^2}
 \le\norm{w_E-w_D}_2.
\end{equation}
Testing the torsion equation of $w_E$ with $w_D$ gives
\[
 \int_E\nabla w_E\cdot\nabla w_D\,dx=\int_Ew_D\,dx=T(D),
\]
hence
\[
 \norm{\nabla(w_E-w_D)}_2^2=T(E)-T(D).
\]
Poincar\'e's inequality on $E$ then yields
\[\norm{w_E-w_D}_2^2\le \dfrac{ \norm{\nabla(w_E-w_D)}_2^2}{\lambda_1(E)}=\dfrac{ T(E)-T(D)}{\lambda_1(E)}\]
which together with \eqref{restorsion} proves \eqref{eq:nested-resolvents}.
\end{proof}

\begin{lemma}\label{lem:remest}
    There exist $\delta_0(n)>0$ and $C_n>0$ such that, if $\delta\le\delta_0$, then 
    \begin{equation}
        \label{rq:rem1}
        \norm{\mathcal{R}_\Omega-\mathcal{R}_{U}}_{L^\infty\to L^2} \le C_n \delta^{\frac{1}{4}},
    \end{equation}
    and
        \begin{equation}
        \label{rq:rem2}
        \norm{\mathcal{R}_{B_{r}}-\mathcal{R}_{B_{r_0}}}_{L^\infty\to L^2}\le C_n \delta^{1/2}.
    \end{equation}
\end{lemma}
\begin{proof}
    By the Kohler--Jobin inequality applied to $U$, we have
    \[\lambda_1(U)T(U)^{\frac{2}{n+2}}\ge \lambda_1(B_{r_0})=\lambda_1(\Omega)-\delta.\]
    While, by \eqref{eq:UOmegaeigcolose}, if $\delta$ is sufficiently small, we have
    \[\lambda_1(U)\le\lambda_1(\Omega)+C_n\delta^{1/2}.\]
    Hence, combining the two estimates, we get 
    \[\begin{split}
    T(\Omega)-T(U)=&\,1-T(U)\le \dfrac{n+2}{2}\left(1-T(U)^\frac{2}{n+2}\right)\\[5pt]
    &\le \dfrac{n+2}{2}\left(1-\dfrac{\lambda_1(\Omega)-\delta}{\lambda_1(\Omega)+C_n\delta^{1/2}}\right)\le C_n \delta^{1/2}.\end{split}\]
    So that, by \autoref{lem:nested-resolvents}, 
    \[\norm{\mathcal{R}_\Omega-\mathcal{R}_{U}}_{L^\infty\to L^2}^2 \le\dfrac{T(\Omega)-T(U)}{\lambda_1(\Omega)}\le C_n \delta^{1/2},\]
    which is exactly \eqref{rq:rem1}.\medskip

    To prove \eqref{rq:rem2} we notice that by \eqref{eq:r-close} 
    \[0\le r_0-r\le C_n\delta,\]
    so that 
    \[T(B_{r_0})-T(B_r)=1-\left(\dfrac{r}{r_0}\right)^{n+2}\le C_n \delta.\]
    Hence, applying \autoref{lem:nested-resolvents} once again, we get
    \[\norm{\mathcal{R}_{B_{r}}-\mathcal{R}_{B_{r_0}}}_{L^\infty\to L^2}^2\le \dfrac{T(B_{r_0})-T(B_r)}{\lambda_1(B_{r_0})}\le C_n \delta.\qedhere\]
\end{proof}

We can finally prove \autoref{thm:intro-KJ}.

\begin{proof}[proof of \autoref{thm:intro-KJ}]
    Let $\Omega$ be a bounded open set with $T(\Omega)=1$. Then by \autoref{prop:torsion-core-gap} and \autoref{lem:remest}, there exist dimensional constants $C_n,\delta_0>0$ and a point $x_\Omega\in\R^n$ such that, if the Kohler--Jobin deficit of $\Omega$, $\delta$, is smaller than $\delta_0$, then 
    \begin{equation}\label{kjqsmall}\norm{\res_\Omega-\res_{B_{r_0}(x_\Omega)}}_{L^\infty\to L^2}^4\le C_n \delta.\end{equation}

    We now observe that, up to increasing the constant $C_n$, estimate \eqref{kjqsmall} holds for every $\delta$. Indeed, for every bounded open set $D$, arguing as in the proof of \autoref{lem:nested-resolvents}, we have
\begin{equation}\notag
 \norm{\res_D}_{L^\infty\to L^2}
 \le \left(\frac{T(D)}{\lambda_1(D)}\right)^{1/2}.
\end{equation}
 Hence
\[
 \norm{\res_\Omega-\res_{B_{r_0}(x_\Omega)}}_{L^\infty\to L^2}^4
 \le \left(\dfrac{1}{\lambda_1(\Omega)^{1/2}}+\dfrac{1}{\lambda_1({B_{r_0})^{1/2}}}\right)^4\le \dfrac{16}{\lambda_1(B_{r_0})^2}.
\]
If $\delta>\delta_0$, then 
\[\norm{\res_\Omega-\res_{B_{r_0}(x_\Omega)}}_{L^\infty\to L^2}^4
 \le  \dfrac{16}{\lambda_1(B_{r_0})^2}\le \dfrac{16}{\delta_0 \lambda_1(B_{r_0})^2}\delta.\]
 That is, exist a dimensional constant $C_n>0$ such that, for every bounded open set $\Omega$ with $T(\Omega)=1$, 
 \begin{equation}\label{KJq-norm}\delta_{\mathrm{KJ}}(\Omega)\ge C_n\inf_{x\in\R^n}\norm{\res_\Omega-\res_{B_{r_0}(x)}}_{L^\infty\to L^2}^4.\end{equation}

 Finally, let $\Omega$ be any bounded open set and let $B^T=B^{T(\Omega)}$ be the open ball having the same torsion as $\Omega$ and let \[\hat{\Omega}=T(\Omega)^{-\frac{1}{n+2}}\Omega= \alpha \Omega \] 
so that $T(\,\hat{\Omega}\,)=1$ and $\alpha B^{T}= B_{r_0}$. Then, by \eqref{KJq-norm}, 
\begin{equation}\label{2tolast}\delta_{\mathrm{KJ}}(\,\hat\Omega\,)\ge c_n \inf_{x\in\R^n}\norm{\res_{\hat\Omega}-\res_{B_{r_0}(x)}}_{L^\infty\to L^2}^4=c_n\inf_{x\in\R^n} \norm{\res_{\alpha\Omega}-\res_{\alpha B^{T}(x)}}_{L^\infty\to L^2}^4. \end{equation}
As the Kohler--Jobin deficit is a scale-invariant quantity, we have the assertion as soon as we determine the scaling property of the resolvent norm.\medskip

Let $D\subset\R^n$ be a bounded open set and let $t\ge0$. For every $f\in L^\infty(tD)$ let $f_t\,\colon x\in D\mapsto f(tx)$, then we have 
\[(\res_{tD}f)(x)=t^2(\res_D f_t)(t^{-1}x).\]
Hence, as $\norm{f}_{\infty}=\norm{f_t}_\infty$, we have 
\[\norm{\res_{tD}}_{L^\infty\to L^2}^2=t^{n+4}\norm{\res_D}_{L^\infty\to L^2}^2.\]

Turning back to \eqref{2tolast} we finally have
\[\delta_{\mathrm{KJ}}(\Omega)=\delta_{\mathrm{KJ}}(\,\hat\Omega\,)\ge c_n\inf_{x\in\R^n} \norm{\res_{\alpha\Omega}-\res_{\alpha B^{T}(x)}}_{L^\infty\to L^2}^4=c_n \left(\inf_{x\in\R^n}\dfrac{\norm{\res_\Omega-\res_{B^{T}(x)}}_{L^\infty\to L^2}}{T(\Omega)^\frac{n+4}{2(n+2)}}\right)^4,
\]
which is the assertion.
\end{proof}

\begin{corollary}\label{cor:KJ-disconnected}
    \autoref{thm:intro-KJ} remains valid for every bounded open set $\Omega\subset\R^n$, without assuming connectedness.
\end{corollary}
\begin{proof}
Let $\Omega$ be a bounded open set. Without loss of generality, assume $T(\Omega)=1$ and $\Omega=\Omega_1\cup\Omega_2$, where $\Omega_1$ is a connected component such that $\lambda_1(\Omega)=\lambda_1(\Omega_1)$. Then we have 
    \[\begin{split}\delta_{\mathrm{KJ}}(\Omega)&= T(\Omega)^\frac{2}{n+2}\lambda_1(\Omega)-T(B)^\frac{2}{n+2}\lambda_1(B)\\[10pt]
    &=\left(1+\dfrac{T(\Omega_2)}{T(\Omega_1)}\right)^\frac{2}{n+2}T(\Omega_1)^\frac{2}{n+2}\lambda_1(\Omega_1)-T(B)^\frac{2}{n+2}\lambda_1(B)\\[10pt]
    &\ge\left[\left(1+\dfrac{T(\Omega_2)}{T(\Omega_1)}\right)^\frac{2}{n+2}-1\right]T(B)^\frac{2}{n+2}\lambda_1(B)+\delta_{\mathrm{KJ}}(\Omega_1)\ge C_n T(\Omega_2)+\delta_{\mathrm{KJ}}(\Omega_1).
    \end{split}\]
    On the other hand, letting $B=B^{T(\Omega)}$ and $B'=B^{T(\Omega_1)}$, we have
    \[\begin{split}\inf_{x\in\R^n}\norm{\res_\Omega-\res_{B(x)}}_{L^\infty\to L^2}&\le\inf_{x\in\R^n}\norm{\res_{\Omega_1}-\res_{B'(x)}}_{L^\infty\to L^2}+\norm{\res_{\Omega_2}}_{L^\infty\to L^2}+\norm{\res_{B'}-\res_B}_{L^\infty\to L^2}\\[10pt]&\le\inf_{x\in\R^n}\norm{\res_{\Omega_1}-\res_{B'(x)}}_{L^\infty\to L^2}+\left(\dfrac{T(\Omega_2)}{\lambda_1(\Omega_2)}\right)^{\frac{1}{2}}+\left(\dfrac{T(B)-T(B')}{\lambda_1(B)}\right)^{\frac{1}{2}}\\[10pt]
    &\le C_n\left(\delta_{\mathrm{KJ}}(\Omega_1)^\frac{1}{4}+T(\Omega_2)^\frac{1}{2}\right)\le C_n\left(\delta_{\mathrm{KJ}}(\Omega_1)+T(\Omega_2)\right)^\frac{1}{4}.
    \end{split}\]
    Putting together the two estimates, we obtain the assertion for disconnected sets.
\end{proof}

\subsection{Further consequences}\label{KJ-fc}

The stability result proved in \autoref{thm:intro-KJ} also provides quantitative control of the entire Dirichlet spectrum and its eigenfunctions.
Indeed, we have the following two results.

\begin{proposition}\label{prop:discussion-spectral}
Let $n\ge2$. There exists a positive, dimensional constant $C_n$ such that, for every bounded, open set $\Omega\subset\mathbb R^n$, 

\begin{equation}\label{eq:discussion-L2-resolvent}
\delta_{\mathrm{KJ}}(\Omega)\ge C_n \left(\inf_{x\in\R^n}\dfrac{\norm{\res_\Omega-\res_{B^{T(\Omega)}(x)}}_{L^2\to L^2}}{T(\Omega)^{\frac2{n+2}}}\right)^{2(n+2)}
\end{equation}
Consequently, for every $j\ge1$,
\begin{equation}\label{eq:discussion-spectral-uniform}
\delta_{\mathrm{KJ}}(\Omega)\ge C_n \left(T(\Omega)^{-\frac2{n+2}}
\abs{
\frac1{\lambda_j(\Omega)}
-\frac1{\lambda_j(B^{T(\Omega)})}
}\right)^{2(n+2)}
\end{equation}
\end{proposition}

\begin{proof}
We start by proving\eqref{eq:discussion-L2-resolvent}. As $\Omega$ is bounded, by \autoref{thm:intro-KJ} and by scaling and translation invariance, we may assume, without loss of generality, that
$T(\Omega)=T(B_{r_0})=1$ and
\begin{equation}\notag
\norm{\res_\Omega-\res_{B_{r_0}}}_{L^\infty\to L^2}
\le C_n\delta_{\mathrm{KJ}}(\Omega)^{\frac{1}{4}}.
\end{equation}

Since the Dirichlet resolvents are self-adjoint, for every
$f\in L^2(\mathbb R^n)$ and
$g\in L^\infty(\mathbb R^n)\cap L^2(\mathbb R^n)$,
\[
\begin{split}
\abs{\langle(\res_\Omega-\res_{B_{r_0}})f,g\rangle_{L^2}}
&=\abs{\langle f,(\res_\Omega-\res_{B_{r_0}})g\rangle_{L^2}}\\[5pt]
&\le \|f\|_2
\norm{(\res_\Omega-\res_{B_{r_0}})g}_2\\[5pt]
&\le C_n\delta_{\mathrm{KJ}}(\Omega)^{\frac{1}{4}}
\norm{f}_2\norm{g}_\infty.
\end{split}
\]
Taking the supremum over $g\in L^\infty(\R^n)$ such that $\norm{g}_\infty\le1$, we have
\begin{equation}\label{eq:discussion-duality}
\norm{(\res_\Omega-\res_{B_{r_0}})f}_1
\le C_n\delta_{\mathrm{KJ}}(\Omega)^{\frac{1}{4}}\norm{f}_2.
\end{equation}

We next control the gradients. By standard elliptic estimates, for every bounded open set $D$,
testing the Dirichlet equation with $\res_D f$ gives
\[
\norm{\nabla\res_D f}_2^2
=\langle f,\res_D f\rangle_{L^2}\le \norm{f}_2\norm{\res_D f}_2\le \lambda_1(D)^{-\frac{1}{2}}\,
\norm{f}_2\norm{\nabla\res_D f}_2.
\]
Hence,
\[
\norm{\nabla\res_D f}_2
\le\lambda_1(D)^{-\frac{1}{2}}\,\norm{f}_2.
\]
Since $T(\Omega)=T(B_{r_0})$, the Kohler--Jobin
inequality ensures that $\lambda_1(\Omega)\ge\lambda_1(B_{r_0})$.
Thus, the previous gradient estimate yields
\begin{equation}\label{eq:discussion-gradient-resolvent}
\norm{\nabla(\res_\Omega-\res_{B_{r_0}})f}_2
\le 2\lambda_1(B_{r_0})^{-\frac{1}{2}}\,\norm{f}_2
\le C_n\norm{f}_2.
\end{equation}

Finally, the Nash inequality
\[
\norm{v}_2
\le C_n\norm{v}_1^{\frac2{n+2}}\,
\norm{\nabla v}_2^{\frac n{n+2}}
\]
applied to $v=(\res_\Omega-\res_{B_{r_0}})f$,
together with \eqref{eq:discussion-duality} and
\eqref{eq:discussion-gradient-resolvent}, gives
\[
\begin{aligned}
\norm{(\res_\Omega-\res_{B_{r_0}})f}_2
&\le C_n
\norm{(\res_\Omega-\res_{B_{r_0}})f}_1^{\frac2{n+2}}\,
\norm{\nabla(\res_\Omega-\res_{B_{r_0}})f}_2^{\frac n{n+2}}\le C_n\delta_{\mathrm{KJ}}(\Omega)^{\frac1{2(n+2)}}\,\norm{f}_2.
\end{aligned}
\]
Taking the supremum over $\norm{f}_2=1$ 
proves \eqref{eq:discussion-L2-resolvent}.\medskip

We now deduce \eqref{eq:discussion-spectral-uniform}
from the $L^2\to L^2$ resolvent estimate.
Recall that, for any bounded open set $D\subset\R^n$, the Dirichlet resolvent is a compact,
non-negative and self-adjoint operator on
$L^2(\mathbb R^n)$, with positive eigenvalues
given by the reciprocals of the Dirichlet eigenvalues.
Thus, from Courant min--max principle we have
\[
\frac1{\lambda_j(D)}
=
\sup_{\substack{V\subset L^2(\mathbb R^n)\\
\dim V=j}}\,
\inf_{\substack{v\in V\\\norm{v}_2=1}}
\langle\res_Dv,v\rangle_{L^2}.
\]

For every $v\in L^2(\mathbb R^n)$ with $\norm{v}_2=1$,
the Cauchy--Schwarz inequality gives
\[
\abs{
\langle\res_\Omega v,v\rangle_{L^2}
-\langle\res_{B_{r_0}}v,v\rangle_{L^2}}=
\abs{\langle
(\res_\Omega-\res_{B_{r_0}})v,v
\rangle_{L^2}}\le
\norm{\res_\Omega-\res_{B_{r_0}}}_{L^2\to L^2}.
\]
Hence, by \eqref{eq:discussion-L2-resolvent}, for every $v\in L^2(\R^n)$ with $\norm{v}_2=1$,
\[
\begin{aligned}
    \langle\res_\Omega v,v\rangle_{L^2}&\le \langle\res_{B_{r_0}}v,v\rangle_{L^2}+C_n \delta_{\mathrm{KJ}}(\Omega)^\frac{1}{2(n+2)},\\[5pt]
    \langle\res_{B_{r_0}} v,v\rangle_{L^2}&\le \langle\res_\Omega v,v\rangle_{L^2}+C_n \delta_{\mathrm{KJ}}(\Omega)^\frac{1}{2(n+2)}.
\end{aligned}
\]
Taking first the infimum over any $j$-dimensional
subspace and then the supremum over all such
subspaces then yields
\[\begin{aligned}
    \frac1{\lambda_j(\Omega)}
&\le
\frac1{\lambda_j(B_{r_0})}
+C_n \delta_{\mathrm{KJ}}(\Omega)^\frac{1}{2(n+2)},\\[10pt]
\frac1{\lambda_j(B_{r_0})}
&\le
\frac1{\lambda_j(\Omega)}
+C_n \delta_{\mathrm{KJ}}(\Omega)^\frac{1}{2(n+2)}.
\end{aligned}\]
That is \eqref{eq:discussion-spectral-uniform}.
\end{proof}

Letting $\mathcal{E}_j(x)=\mathcal E_j(B^{T(\Omega)}(x))\subset H^1_0(B^{T(\Omega)}(x))$ denote the
eigenspace corresponding to the eigenvalue $\lambda_j(B^{T(\Omega)})$, we have the following.
\begin{proposition}
Let $n\ge2$. For every $j\ge1$, there exists a positive, dimensional constant $C_{n,j}$ such that, for every bounded, open set $\Omega\subset\mathbb R^n$ and every $L^2$-normalized $j$-th Dirichlet eigenfunction $u_{\Omega,j}$,
\begin{equation}\notag
\begin{aligned}
\delta_{\mathrm{KJ}}(\Omega)&\ge C_{n,j}\left(T(\Omega)^{-\frac2{n+2}}
\abs{\dfrac{1}{\lambda_j(\Omega)}
-\dfrac{1}{\lambda_j(B^{T(\Omega)})}}\right)^4
,\\[10pt]
\delta_{\mathrm{KJ}}(\Omega)&\ge C_{n,j} \left(\inf_{\substack{
v\in\mathcal E_j(x_\Omega)\\
\norm{v}_2=1}}
\|u_{\Omega,j}-v\|_{L^2(\mathbb R^n)}\right)^4,
\end{aligned}
\end{equation}
for some $x_\Omega\in\R^n$ independent of $j$.
\end{proposition}
\begin{proof}
Fix $j\ge1$. By the Kohler--Jobin inequality 
$\lambda_1(\Omega)\ge\lambda_1(B_{r_0})$, so that
\[
\abs{
\frac1{\lambda_j(\Omega)}
-\frac1{\lambda_j(B_{r_0})}
}
\le\frac2{\lambda_1(B_{r_0})}.
\]
similarly, the $L^2$-distance between two
normalized functions is at most $2$. Hence it is
enough to establish both estimates when
$\delta_{\mathrm{KJ}}(\Omega)$ is sufficiently small,
depending only on $n$ and $j$. As in the proof of \autoref{prop:discussion-spectral}, we assume, without loss of generality, that
$T(\Omega)=T(B_{r_0})=1$ and
\begin{equation}\label{eq:assumption-prop-dist-eigfunction}
\norm{\res_\Omega-\res_{B_{r_0}}}_{L^\infty\to L^2}
\le C_n\delta_{\mathrm{KJ}}(\Omega)^{\frac{1}{4}}.
\end{equation}

We start by estimating the reciprocal of the eigenvalues.
Let $\{\varphi_k\}_{k\ge1}$ be an $L^2$-orthonormal
basis of Dirichlet eigenfunctions of $B_{r_0}$,
extended by zero to $\mathbb R^n$, and set
\[
V=\operatorname{span}\{\varphi_1,\ldots,\varphi_j\}.
\]
Let $v\in V$, then there exists $c_1,\dots,c_j\in \R$  such that
$v=\sum_{i=1}^j c_i\varphi_i$. Thus, by the Cauchy--Schwarz inequality, we have
\[
\begin{aligned}
\norm{v}_\infty
&\le
\left(\sum_{i=1}^j\norm{\varphi_i}_\infty^2\right)^{1/2}
\left(\sum_{i=1}^j\abs{c_i}^2\right)^{1/2}\le C_{n,j}\norm{v}_2,
\end{aligned}
\]
where we have used that the $L^\infty$ norm of $\varphi_i$ can be seen as constants only depending on $n$ and $j$.\medskip

Thus, for every $v\in V$ such that $\norm{v}_2=1$, we can deduce
\[
\begin{split}
\langle\res_\Omega v,v\rangle_{L^2}
&\ge
\langle\res_{B_{r_0}}v,v\rangle_{L^2}
-\norm{(\res_\Omega-\res_{B_{r_0}})v}_2\ge\\
&\ge
\langle\res_{B_{r_0}}v,v\rangle_{L^2}
-C_n\delta_{\mathrm{KJ}}(\Omega)^{\frac{1}{4}}\norm{v}_\infty\ge
\langle\res_{B_{r_0}}v,v\rangle_{L^2}
-C_{n,j}\delta_{\mathrm{KJ}}(\Omega)^{\frac{1}{4}}.
\end{split}
\]
By the choice of $V$, we have
\[
\inf_{\substack{v\in V\\\|v\|_2=1}}
\langle\res_{B_{r_0}}v,v\rangle_{L^2}
=\frac1{\lambda_j(B_{r_0})},
\]
hence, the min--max principle for $\res_\Omega$ yields
\[
\frac1{\lambda_j(\Omega)}
\ge
\frac1{\lambda_j(B_{r_0})}
-C_{n,j}\delta_{\mathrm{KJ}}(\Omega)^{\frac{1}{4}}.
\]
For sufficiently small deficit, this implies
\[
\lambda_j(\Omega)\le2\lambda_j(B_{r_0}).
\]

Therefore, we can apply the universal eigenfunction
estimate in \cite[Theorem~1.6, equation~(1.8)]
{vanDenBergHempelVoigt2015}
\[
\norm{u_{\Omega,i}}_\infty
\le C_n\lambda_i(\Omega)^{n/4}\norm{u_{\Omega,i}}_2,
\]
to obtain
\begin{equation}\label{linfinty up bound}
   \norm{u_{\Omega,i}}_\infty
\le C_n(2\lambda_j(B_{r_0}))^{n/4}
\le C_{n,j}. 
\end{equation}
As a consequence, since the eigenfunctions of $ \Omega$ are uniformly bounded by constants depending only on $n$ and $j$, we can repeat the same argument interchanging $\Omega$ and $B_{r_0}$. Thus we obtain
\[
\frac1{\lambda_j(B_{r_0})}
\ge
\frac1{\lambda_j(\Omega)}
-C_{n,j}\delta_{\mathrm{KJ}}(\Omega)^{\frac{1}{4}}.
\]
Combining the two inequalities gives
\begin{equation}\label{eq:discussion-reciprocal-small}
\left|
\frac1{\lambda_j(\Omega)}
-\frac1{\lambda_j(B_{r_0})}
\right|
\le C_{n,j}\delta_{\mathrm{KJ}}(\Omega)^{\frac{1}{4}}.
\end{equation}

Following the argument underlying
\cite[proofs of Theorems~5.3 and~5.4]{AKN25}, we now prove the quantitative control of eigenfunctions.
Since
\[
\res_\Omega u_{\Omega,j}
=\frac1{\lambda_j(\Omega)}u_{\Omega,j},
\]
using \eqref{linfinty up bound}, \eqref{eq:discussion-reciprocal-small}  and \eqref{eq:assumption-prop-dist-eigfunction} we obtain

\begin{equation}\label{uOmBeig}
\norm{
\res_{B_{r_0}}u_{\Omega,j}
-\frac1{\lambda_j(B_{r_0})}u_{\Omega,j}
}_2\le
\norm{(\res_{B_{r_0}}-\res_\Omega)u_{\Omega,j}}_2
+\abs{
\frac1{\lambda_j(\Omega)}
-\frac1{\lambda_j(B_{r_0})}
}
\le C_{n,j}\delta_{\mathrm{KJ}}(\Omega)^{\frac{1}{4}}.
\end{equation}

Let $P_j$ denote the orthogonal projection of
$L^2(\mathbb R^n)$ onto $\mathcal E_j$.
Using the notations introduced above, we have
\[\begin{aligned}
    &f= f\mathbf1_{\R^n\setminus B_{r_0}}+\sum_{k\ge1} \langle f,\varphi_k\rangle_{L^2}\,\varphi_k,\\[10pt]
    &\res_{B_{r_0}} f = \sum_{k\ge1} \dfrac{\langle f,\varphi_k\rangle_{L^2}}{\lambda_k (B_{r_0}) }\,\varphi_k,\\[10pt]
    &P_jf = \sum_{\substack{k\ge1\\ \lambda_k(B_{r_0})=\lambda_j(B_{r_0})}} \langle f,\varphi_k\rangle_{L^2}\,\varphi_k.
\end{aligned}\]
Hence, we obtain
\begin{equation}\label{aa}
\begin{split}
\norm{
\res_{B_{r_0}}f-\frac1{\lambda_j(B_{r_0})}f
}_2^2
&=
\sum_{\substack{k\ge1\\
\lambda_k(B_{r_0})\neq\lambda_j(B_{r_0})}}
\abs{
\frac1{\lambda_k(B_{r_0})}
-\frac1{\lambda_j(B_{r_0})}
}^2
\abs{\langle f,\varphi_k\rangle_{L^2}}^2+
\frac1{\lambda_j(B_{r_0})^2}
\norm{f}_{L^2(\mathbb R^n\setminus B_{r_0})}^2\\[10pt]
&\ge c_{n,j}^2\norm{f-P_jf}_2^2.
\end{split}
\end{equation}
where we have used that, since $\res_{B_{r_0}}$ is compact and self-adjoint,
$1/\lambda_j(B_{r_0})$ is isolated from the rest of its spectrum.
Therefore, we have that 
\[
c_{n,j}:= 
\operatorname{dist}\left(
\lambda_j(B_{r_0})^{-1},
\sigma(\res_{B_{r_0}})\setminus
\Set{\lambda_j(B_{r_0})^{-1}}
\right)>0.
\]
can be seen as a positive constant depending only on $n$ and $j$.\medskip

Applying \eqref{aa} to $u_{\Omega,j}$ and using \eqref{uOmBeig},
we obtain
\begin{equation}\notag
\norm{u_{\Omega,j}-P_ju_{\Omega,j}}_2
\le C_{n,j}\delta_{\mathrm{KJ}}(\Omega)^{\frac{1}{4}}.
\end{equation}
Finally, since $u_{\Omega,j}$ is normalized and $P_j$
is an orthogonal projection,
\[
1\ge\norm{P_ju_{\Omega,j}}_2^2
=1-\norm{u_{\Omega,j}-P_ju_{\Omega,j}}_2^2
\ge1-C_{n,j}\delta_{\mathrm{KJ}}(\Omega)^{1/2}.
\]
Thus, for sufficiently
small deficit, $P_ju_{\Omega,j}\neq0$ and we can define
\[
v=\frac{P_ju_{\Omega,j}}
{\|P_ju_{\Omega,j}\|_2}
\in\mathcal E_j.
\]
By construction, $\norm{v}_2=1$, and, 
by orthogonality,
\[
\begin{split}
\norm{u_{\Omega,j}-v}_2^2
&=2\bigl(1-\norm{ P_j u_{\Omega,j}}_2\bigr)\le2\bigl(1-\norm{P_ju_{\Omega,j}}_2^2\bigr)\\[5pt]
&=2\norm{u_{\Omega,j}-P_ju_{\Omega,j}}_2^2\le C_{n,j}\delta_{\mathrm{KJ}}(\Omega)^{\frac1{2}}
\end{split}
\]
which concludes the proof.
\end{proof}

Finally, as noted in the introduction, the resolvent stability estimate of \autoref{thm:intro-KJ} also provides information on the
convergence of the domains themselves. The natural framework is the
global $\gamma$-distance introduced in \cite[Section~2.5]{BBV}.

\begin{definition}\label{defi:discussion-gamma}
Let $\Omega_1,\Omega_2\subset\mathbb R^n$ be two bounded open sets, and let $w_{\Omega_1}$ and $w_{\Omega_2}$, denote the torsion function of $\Omega_1$ and $\Omega_2$ respectively. We define the $\gamma$-distance between $\Omega_1$ and $\omega_2$ as
\[
d_\gamma(\Omega_1,\Omega_2)
:=\|w_{\Omega_1}-w_{\Omega_2}\|_{L^1(\mathbb R^n)},
\]
where, as usual, the torsion functions are extended by zero outside their
respective domains. Thus, a sequence $(\Omega_k)$ is said to
$\gamma$-converge to $\Omega$ if
\[
d_\gamma(\Omega_k,\Omega)\longrightarrow0.
\]
\end{definition}

This is the restriction to bounded open sets of the distance
defined in \cite[Section~2.5]{BBV} for capacitary measures
of finite torsion. In particular, we remark that, if $(\Omega_k)$ $\gamma$-convergence to $\Omega$ then
\[
T(\Omega_k)\longrightarrow T(\Omega)
\]
and, by \cite[Remark~2.8]{BBV},
\begin{equation}\notag
\|\res_{\Omega_k}-\res_\Omega\|_{L^2\to L^2}
\longrightarrow0.
\end{equation}

We can now prove \autoref{thm-KJ-gamma}
\begin{proof}[proof of \autoref{thm-KJ-gamma}]
Let $B=B^{T(\Omega)}(x)$ be any translate of the ball with the
same torsional rigidity as $\Omega$. Since $w_\Omega$ and $w_B$
are non-negative and have the same integral, we have
\[
\begin{aligned}
d_\gamma(\Omega,B)
&=\int_{\mathbb R^n}\abs{w_\Omega-w_B}\,dx=2\int_{\mathbb R^n}(w_B-w_\Omega)_+\,dx=2\int_B(w_B-w_\Omega)_+\,dx\le 2\abs{B}^{\frac1{2}}\norm{w_\Omega-w_B}_2.
\end{aligned}
\]
Moreover,  the definition of the resolvent norm gives
\[
\norm{w_\Omega-w_B}_2=\norm{\left(\res_\Omega-\res_B\right)1}_2
\le\norm{\res_\Omega-\res_B}_{L^\infty\to L^2},
\]
so that
\[
d_\gamma(\Omega,B)
\le 2\abs{B}^{\frac1{2}}
\norm{\res_\Omega-\res_B}_{L^\infty\to L^2}.
\]
Since
\[
\abs{B}
=\abs{B_{r_0}}T(\Omega)^{\frac{n}{n+2}},
\]
 \autoref{thm:intro-KJ} ensures
\[
\begin{aligned}
\inf_{x\in\mathbb R^n}
d_\gamma\left(\Omega,B^{T(\Omega)}(x)\right)
&\le C_n T(\Omega)^{\frac{n}{2(n+2)}}
T(\Omega)^{\frac{n+4}{2(n+2)}}
\delta_{\mathrm{KJ}}(\Omega)^{\frac{1}{4}}=C_nT(\Omega)\delta_{\mathrm{KJ}}(\Omega)^{\frac{1}{4}},
\end{aligned}
\]
which proves \eqref{eq:discussion-KJ-gamma-quantitative:intro}.\medskip

Let $(\Omega_k)$ be a sequence of bounded open sets such that $\delta_{\mathrm{KJ}}(\Omega_k)\to0$.
Choose translations $x_k$ so that
\[d_\gamma\left(\Omega_k,B^{T(\Omega_k)}(x_k)\right)\le C_n T(\Omega_k)\delta_{\mathrm{KJ}}(\Omega_k)^{\frac{1}{4}},\]
and set
\[
E_k=\frac{\Omega_k-x_k}{r_k},
\]
for 
\[r_k=\left(\dfrac{T(\Omega_k)}{T(B_1)}\right)^\frac1{n+2}.\]
Recalling that
\[
w_{rE}(x)=r^2w_E(x/r),
\]
we have that, the $\gamma$-distance satisfies
\[
d_\gamma(rE,rF)=r^{n+2}d_\gamma(E,F).
\]
Consequently,
\[
\begin{aligned}
d_\gamma(E_k,B_1)
&=\frac{T(B_1)}{T(\Omega_k)}
d_\gamma\left(\Omega_k,B^{T(\Omega_k)}(x_k)\right)\le C_n\delta_{\mathrm{KJ}}(\Omega_k)^{\frac{1}{4}}
\longrightarrow0.
\end{aligned}
\]
That is $E_k$ $\gamma-$converges to $B_1$.\medskip

Conversely, suppose that $E_k\xrightarrow{\gamma}B_1$. Then, as stated,
by \cite[Remark~2.8]{BBV}, $\res_{E_k}$
converge in $L^2$ operator norm to $\res_{B_1}$, hence in particular
\[
\lambda_1(E_k)\longrightarrow\lambda_1(B_1).
\]
Since, by definition, $T(E_k)=T(B_1)$, we have 
\[
\delta_{\mathrm{KJ}}(\Omega_k)
=\delta_{\mathrm{KJ}}(E_k)
=T(B_1)^{\frac{2}{n+2}}(\lambda_1(E_k)-\lambda_1(B_1))
\longrightarrow0.
\]
This proves the assertion.
\end{proof}

\section{On the different nature of the two stability results}\label{sec:counterexamples}

Although the proofs of \autoref{thm:intro-Chiti} and \autoref{thm:intro-KJ} follow a common strategy, the quantitative inequalities proved involve two different notions of
asymmetry. Namely, the Chiti deficit controls the distance between the normalized
first Dirichlet eigenfunction of $\Omega$ and that of the ball $B^\lambda=B^{\lambda_1(\Omega)}$
with $\lambda_1(B^\lambda)=\lambda_1(\Omega)$, whereas the
Kohler--Jobin deficit controls the distance between the Dirichlet
resolvents of $\Omega$ and a ball $B^T=B^{T(\Omega)}$ having the same
torsional rigidity. 
The latter provides a stronger form of stability,
since it compares the solutions of the Dirichlet problem for every
bounded right-hand side.\medskip

Such a resolvent estimate cannot hold for the Chiti deficit.
Indeed, consider $\Omega$ the union of two disjoint balls of radii $R$ and $r$ with $R>r$ such that $\operatorname{dist}(B_R,B_r)>4R$.
Since $\lambda_1(\Omega)=\lambda_1(B_R)$, a non-negative first
Dirichlet eigenfunction of $\Omega$ coincides
with the first eigenfunction of $B_R$, extended by zero
on the other ball. Consequently, the Chiti deficit vanishes.

On the other hand, every translate $B_R(x)$ is disjoint from
at least one of the two balls, say $F$. Taking
$f=\mathbf 1_F$, we have
\[
\mathcal R_\Omega f=w_F,
\qquad
\mathcal R_{B_R(x)}f=0.
\]
It follows that
\begin{equation}\notag
\inf_{x\in\mathbb R^n}
\norm{\mathcal R_\Omega-\mathcal R_{B_R(x)}}_{L^\infty\to L^2}
\ge \norm{w_{B_r}}_{L^2(\mathbb R^n)}>0.
\end{equation}
Thus, even when the Chiti deficit is zero, the resolvent
may remain at a positive distance from that of the
comparison ball.\medskip

The same obstruction persists for connected domains. Indeed, consider $\Omega$ composed of two disjoint balls of fixed radii \(R>r>0\), joined by a thin tube whose width tends to zero. By the standard spectral convergence for shrinking handles, the first eigenvalue of $\omega$ converges to \(\lambda_1(B_R)\), and the corresponding normalized eigenfunction converges to the first eigenfunction of the larger ball, extended by zero; see also \cite{AFT14}. Consequently, the Chiti deficit tends to zero, whereas the resolvent distance from any translated comparison ball remains bounded away from zero, since the smaller ball is still detected by the resolvent.

\begin{figure}[H]
\centering
\begin{tikzpicture}[scale=1.15,
    line cap=round, line join=round]

    \path[
        fill=blue!15,
        draw=blue!65!black,
        line width=0.85pt
    ]
    (-0.653,0.08)
    arc[start angle=4.59,end angle=355.41,radius=1]
    -- (0.774,-0.08)
    arc[start angle=185.89,end angle=534.11,radius=0.78]
    -- cycle;

    \node at (-1.65,0) {$B_R$};
    \node at (1.55,0) {$B_r$};

    \draw[<->, thin]
        (0.10,-0.08) -- (0.10,0.08);
    \node[above] at (0.10,0.18) {$\varepsilon$};

\end{tikzpicture}
\caption{Two balls of fixed radii $R>r$ connected
by a thin tube of width $\varepsilon\to0$.
The first normalized Dirichlet eigenfunction
concentrates on the larger ball, whereas
the smaller component remains detectable
by the resolvent.}
\label{fig:discussion-dumbbell}
\end{figure}
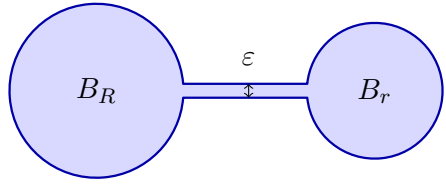

A natural question is whether the asymmetries appearing
in the preceding quantitative estimates can be replaced
by the more classical Fraenkel asymmetry. This geometric
notion of distance from a ball plays a central role in
quantitative inequalities, notably in the quantitative
Faber--Krahn and Saint--Venant inequalities; see
\cite{BDPV15}. However, unlike these classical inequalities,
neither Chiti's nor Kohler--Jobin's inequality is formulated
under a prescribed-volume constraint. Indeed, the comparison
ball is determined by the first Dirichlet eigenvalue in the
former case and by the torsional rigidity in the latter.
It is therefore more natural to measure the distance
between the corresponding PDE quantities rather than
the distance between the domains in measure. In fact, we can show that
neither deficit controls the Fraenkel asymmetry.\medskip

For Chiti's inequality, this already follows from the
dumbbell construction discussed above
(see \autoref{fig:discussion-dumbbell}). Indeed, as the connecting
tube shrinks, the Chiti deficit tends to zero, while the
Fraenkel asymmetry remains bounded away from zero.\medskip

For the Kohler--Jobin inequality, the same construction
does not provide a counterexample, since the smaller ball
retains a positive contribution to the torsional rigidity.
Instead, the obstruction is related to the fact that
$\gamma$-convergence does not necessarily imply convergence
in measure. Indeed, our quantitative Kohler--Jobin
inequality guarantees that, up to translations and
dilations, domains with vanishing deficit $\gamma$-converge
to a ball. Nevertheless, such convergence does not
force the volume outside the limiting ball to vanish.
This phenomenon is well known in shape optimization;
see, for instance, \cite[Introduction]{BucurDurusOudet}
and \cite[Section~2.5]{BBV}.\medskip

The following construction makes this distinction
explicit, showing that the Kohler--Jobin deficit may
vanish even when the Fraenkel asymmetry converges to
any prescribed value in $(0,2)$.

\begin{proposition}
Let $n\ge2$ and $\alpha\in(0,2)$. There exist $R>1$ and
a sequence of bounded, connected Lipschitz domains
$\Omega_m\subset\mathbb R^n$ such that
\[
B_1\subset\Omega_m\subset B_R,
\qquad
d_\gamma(\Omega_m,B_1)\longrightarrow0,
\qquad
\mathcal A(\Omega_m)\longrightarrow\alpha,
\]
where $d_\gamma$ is the distance introduced in
\autoref{defi:discussion-gamma}.
In particular,
\[
\delta_{\mathrm{KJ}}(\Omega_m)\longrightarrow0.
\]
\end{proposition}

\begin{proof}
We construct the set $\Omega_m$ by attaching an increasing number of increasingly
thin radial tentacles to $B_1$.\medskip

Fix $R>1$ and $\vartheta\in(0,1)$ that will be chosen later on in the proof, and let $\{\varepsilon_m\}$ be a vanishing sequence of positive real numbers. For every $m$, consider a finite partition of $\mathbb S^{n-1}$,
up to a set of $\mathcal{H}^{n-1}$-measure zero, into Lipschitz sets
$Q_{m,1},\ldots,Q_{m,N_m}$ satisfying
\[
\max_{1\le j\le N_m}\operatorname{diam}(Q_{m,j})
\le\varepsilon_m.
\]
For every $m$ and every $1\le j\le N_m$ let $U_{m,j}$ be a Lipschitz open such that
$U_{m,j}\subset\subset Q_{m,j}$, and
\[
\mathcal{H}^{n-1}(U_{m,j})=\vartheta\,\mathcal{H}^{n-1}(Q_{m,j}).
\]
In particular, by construction, the closures of the sets $U_{m,j}$ are
pairwise disjoint.\medskip

Set
\[
E_m=\bigcup_{j=1}^{N_m}U_{m,j}.
\]
Writing points in polar coordinates as $x=r\xi$, where
$r=|x|>0$ and $\xi=x/|x|\in\mathbb S^{n-1}$, define
\[
\Omega_m
=
B_1\cup
\bigl\{r\xi:\ 1-\eta<r<R,\ \xi\in E_m\bigr\},
\]
where $\eta\in(0,1)$ is fixed.

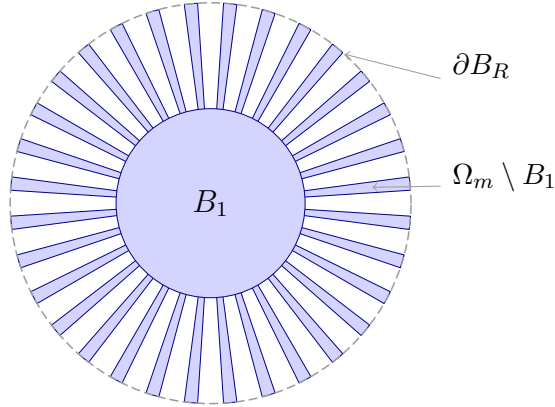
\begin{figure}[htbp]
\centering
\begin{tikzpicture}[
    scale=1.25,
    line join=round,
    line cap=round
]
    \def\OuterRadius{2.12}
    \def\NumberOfCells{32}
    \def\AngularFraction{0.34}

    \foreach \k in {0,...,31} {
        \pgfmathsetmacro{\leftAngle}
            {360*(\k+0.5-\AngularFraction/2)/\NumberOfCells}
        \pgfmathsetmacro{\rightAngle}
            {360*(\k+0.5+\AngularFraction/2)/\NumberOfCells}

        \path[
            fill=blue!17,
            draw=blue!55!black,
            line width=0.32pt
        ]
            (\leftAngle:0.94)
            -- (\leftAngle:\OuterRadius)
            arc[
                start angle=\leftAngle,
                end angle=\rightAngle,
                radius=\OuterRadius
            ]
            -- (\rightAngle:0.94)
            -- cycle;
    }

    \filldraw[
        fill=blue!17,
        draw=blue!55!black,
        line width=0.0pt
    ]
        (0,0) circle[radius=1];

    \draw[
        gray!75,
        densely dashed,
        line width=0.6pt
    ]
        (0,0) circle[radius=\OuterRadius];

    \node[font=\large] at (0,0) {$B_1$};

    \node[anchor=west] at (2.45,1.45) {$\partial B_R$};
    \draw[gray!75, thin, ->]
        (2.43,1.32) -- (48:\OuterRadius);

    \node[anchor=west] at (2.45,0.28)
        {$\Omega_m\setminus B_1$};
    \draw[gray!75, thin, ->]
        (2.43,0.18) -- (5.625:1.72);

\end{tikzpicture}
\caption{Schematic representation of a ball with thin radial
tentacles. As the number of appendages increases and their
angular width decreases, the domains converge to $B_1$ in the
$\gamma$-sense while retaining a positive additional volume.}
\label{fig:discussion-tentacles}
\end{figure}

$\Omega_m$ is a bounded connected Lipschitz domain,
with $B_1\subset\Omega_m\subset B_R$.

Since
\[
\Omega_m\setminus B_1
=\{r\xi:\ 1<r<R,\ \xi\in E_m\},
\]
we obtain
\[
\begin{aligned}
|\Omega_m\setminus B_1|
&=\int_1^R r^{n-1}\,dr
  \int_{E_m}d\mathcal{H}^{n-1}\\[10pt]
&=\frac{R^n-1}{n}\,
  \vartheta\,\mathcal{H}^{n-1}(\mathbb S^{n-1})=\vartheta |B_1|(R^n-1).
\end{aligned}
\]
Consequently, for every $m$
\begin{equation}\label{eq:discussion-sun-volume}
|\Omega_m|
=|B_1|\bigl(1+\vartheta(R^n-1)\bigr).
\end{equation}

We point out that the measure of $\Omega_m\setminus B_1$
is positive and independent of $m$. Moreover, for all $m$
\[\Omega_m\subset\subset B_{R+1}.\]
The proof is divided into two steps. We first compute the limit of the Fraenkel Asymmetry of $\Omega_m$ and then show that $\Omega_m$ $\gamma$-converges to $B_1$.\medskip

To compute the Fraenkel asymmetry, we first show that
\begin{equation}\label{eq:discussion-sun-weakstar-sphere}
\mathbf1_{E_m}\stackrel{*}{\rightharpoonup}\vartheta
\qquad\text{in }L^\infty(\mathbb S^{n-1}).
\end{equation}
Let $g\in C(\mathbb S^{n-1})$, and for every $m$ and $1\le j\le N_m$ fix points
$\xi_{m,j}\in Q_{m,j}$.
Then, as
$\mathcal{H}^{n-1}(U_{m,j})=\vartheta\mathcal{H}^{n-1}(Q_{m,j})$,
we have
\[
\begin{split}
\left|
\int_{\mathbb S^{n-1}}
(\mathbf1_{E_m}-\vartheta)g\,d\mathcal{H}^{n-1}
\right|&=
\left|
\sum_{j=1}^{N_m}
\int_{Q_{m,j}}
(\mathbf1_{U_{m,j}}-\vartheta)
\bigl(g(\xi)-g(\xi_{m,j})\bigr)
\,d\mathcal{H}^{n-1}
\right|\\[10pt]
&\le
\mathcal{H}^{n-1}(\mathbb S^{n-1})
\sup_{\substack{\xi,\xi'\in\mathbb S^{n-1}\\
|\xi-\xi'|\le\varepsilon_m}}
|g(\xi)-g(\xi')|
\longrightarrow0.
\end{split}
\]
By density, we have the claim.\medskip

The convergence \eqref{eq:discussion-sun-weakstar-sphere} implies
\begin{equation}\label{eq:discussion-sun-weakstar-domain}
\mathbf1_{\Omega_m}\stackrel{*}{\rightharpoonup}
\mathbf1_{B_1}
+\vartheta\mathbf1_{B_R\setminus B_1}
\qquad\text{in }L^\infty(\mathbb R^n).
\end{equation}
Indeed, for every $\psi\in L^1(\mathbb R^n)$, we have that
\[
\begin{aligned}
\int_{\mathbb R^n}
\bigl(\mathbf1_{\Omega_m}
-(\mathbf1_{B_1}
+\vartheta\mathbf1_{B_R\setminus B_1})\bigr)
\psi\,dx
&=
\int_{\mathbb S^{n-1}}
(\mathbf1_{E_m}(\theta)-\vartheta)
\left(
\int_1^R\psi(r\theta)r^{n-1}\,dr
\right)d\sigma(\theta).
\end{aligned}
\]
Since the function
\[
\theta\longmapsto
\int_1^R\psi(r\theta)r^{n-1}\,dr
\]
belongs to $L^1(\mathbb S^{n-1})$, by
\eqref{eq:discussion-sun-weakstar-sphere}, we obtain
proving \eqref{eq:discussion-sun-weakstar-domain}.\medskip

Let $B_\rho$ be the ball having the same volume
as $\Omega_m$. By \eqref{eq:discussion-sun-volume},
\[
\rho^n=1+\vartheta(R^n-1),
\qquad 1<\rho<R.
\]
In particular, $\rho$ is independent of $m$.
To compute the Fraenkel asymmetry, we need to
maximize the intersection of $\Omega_m$ with
translations of $B_\rho$.

For $a\in\mathbb R^n$, set
\[
I_m(a)=|\Omega_m\cap B_\rho(a)|=
\int_{\mathbb R^n}
\mathbf1_{\Omega_m}\mathbf1_{B_\rho(a)}\,dx.
\]
Since $\mathbf1_{B_\rho(a)}\in L^1(\mathbb R^n)$,
the weak-* convergence
\eqref{eq:discussion-sun-weakstar-domain}
implies, for every fixed $a\in\mathbb R^n$,
\[
I_m(a)\longrightarrow I(a)=
(1-\vartheta)|B_1\cap B_\rho(a)|
+\vartheta|B_R\cap B_\rho(a)|.
\]

We claim that this convergence is uniform with
respect to $a$. First, we notice that
\[
0\le I_m(a)\le |B_\rho|.
\]
so that the functions $I_m$ are
uniformly bounded. Moreover the functions $I_m$ are
equicontinuous, indeed, for any $a,a'\in\mathbb R^n$, we have
\[
\begin{aligned}
|I_m(a)-I_m(a')|
&\le |B_\rho(a)\triangle B_\rho(a')|\le C(n,\rho)|a-a'|.
\end{aligned}
\]
Finally, for every $\abs{a}>R+\rho$
\[
I_m(a)=0=I(a).
\]
Therefore, we restrict attention to the
compact set $\overline B_{R+\rho}$. On such a set,
the functions $I_m$ are uniformly bounded and
equicontinuous. Hence, by the Arzel\`a--Ascoli theorem,

\begin{equation}\label{eq:discussion-sun-uniform-intersection}
 \lim_{m}
\sup_{a\in\mathbb R^n}I_m(a)
=\sup_{a\in\mathbb R^n}I(a). 
\end{equation}

For every $a\in\mathbb R^n$, we have
\[
|B_1\cap B_\rho(a)|\le |B_1|,
\qquad
|B_R\cap B_\rho(a)|\le |B_\rho|,
\]
with equality for $a=0$.
Hence
\[
I(a)\le
(1-\vartheta)|B_1|+\vartheta|B_\rho|=I(0).
\]
Combining this identity with
\eqref{eq:discussion-sun-uniform-intersection},
we finally obtain
\[
\begin{aligned}
\lim_{m\to\infty}
\sup_{a\in\mathbb R^n}
|\Omega_m\cap B_\rho(a)|
&=(1-\vartheta)|B_1|+\vartheta|B_\rho|=|B_1|\bigl(1+\vartheta^2(R^n-1)\bigr),
\end{aligned}
\]
and, by definition, 
\[
\begin{split}
\lim_{m}\mathcal A(\Omega_m)
&=\lim_{m}\left(\inf_{a} \dfrac{\abs{\Omega_m\triangle B_\rho(a)}}{\abs{\Omega_m}}\right)\\[10pt]&=
2\left(
1-\frac{|\Omega_m\cap B_\rho|}
{|\Omega_m|}
\right)=\frac{2\vartheta(1-\vartheta)(R^n-1)}
{1+\vartheta(R^n-1)}.
\end{split}
\]

Finally, choosing any $\theta$ in $(0,1-\alpha/2)$ and letting
\[
R=
\left(
1+\frac{\alpha}
{\vartheta[2(1-\vartheta)-\alpha]}
\right)^{1/n}.
\]
we have $R>1$, and
$$\lim_{m}\mathcal A(\Omega_m)=\alpha,$$ completing the proof.\medskip

We now establish the $\gamma$-convergence. Since, for every $m$,
$B_1\subset\Omega_m\subset B_R$,  we work in the fixed bounded domain $D=B_{R+1}$,
identifying $H_0^1(\Omega_m)$ and $H_0^1(B_1)$ with subspaces of $H_0^1(D)$ by extension by zero. By \cite[Definition~3.5.4 and Proposition~3.5.5]{HP18shape}, it is sufficient to prove that $H_0^1(\Omega_m)$ converges
to $H_0^1(B_1)$ in \emph{the sense of Mosco}, that is:

\begin{enumerate}
\item For every $u\in H_0^1(B_1)$, there exists
a sequence $u_m\in H_0^1(\Omega_m)$ such that
\[
u_m\longrightarrow u
\qquad\text{strongly in }H_0^1(D).
\]

\item If $u_m\in H_0^1(\Omega_m)$ and
\[
u_m\rightharpoonup u
\qquad\text{weakly in }H_0^1(D),
\]
then $u\in H_0^1(B_1)$.
\end{enumerate}

By the inclusion
$B_1\subset\Omega_m$, the first condition is trivially verified taking the constant sequence $u_m=u\in H^1_0(B_1)\subset H^1_0(\Omega_m)$.\medskip

To verify the second condition, we first show that
functions in $H_0^1(\Omega_m)$ cannot retain a
significant amount of $L^2$-mass in $\Omega_m\setminus B_1$ if their Dirichlet energies remain bounded.
More precisely, we claim that for every $u\in H_0^1(\Omega_m)$
\begin{equation}\label{eq:discussion-sun-poincare}
\int_{\Omega_m\setminus B_1}|u|^2\,dx
\le C(n,R)\varepsilon_m^2
\int_{\Omega_m\setminus B_1}|\nabla u|^2\,dx.
\end{equation}

To prove the claim, recall that the diameter of each angular
set $U_{m,j}$ is smaller than the vanishing sequence $\varepsilon_m$. Thus, $U_{m,j}$ is contained in a spherical cap
of radius $\varepsilon_m$. If $m$ is sufficiently large then, the Dirichlet Poincar\'e inequality on such
spherical caps gives that, for every $h\in H_0^1(U_{m,j})$,
\[
\int_{U_{m,j}}|h|^2\,d\mathcal H^{n-1}
\le C_n\varepsilon_m^2
\int_{U_{m,j}}
|\nabla_{\mathbb S^{n-1}}h|^2\,d\mathcal H^{n-1},
\]
with a constant independent of $m$ and $j$.\medskip

Let $u\in C_c^\infty(\Omega_m)$.
For every $r\in(1,R)$, the function
$\xi\mapsto u(r\xi)$ belongs to
$H_0^1(U_{m,j})$ on each angular section, and
\[
\left|
\nabla_{\mathbb S^{n-1}}
\bigl(u(r\,\cdot)\bigr)(\xi)
\right|
\le r|\nabla u(r\xi)|.
\]
Applying the preceding Poincar\'e inequality
on each $U_{m,j}$ and integrating in polar
coordinates, we obtain
\[
\begin{aligned}
\int_{\Omega_m\setminus B_1}|u|^2\,dx
&=\int_1^R r^{n-1}
\sum_{j=1}^{N_m}
\int_{U_{m,j}}|u(r\xi)|^2
\,d\mathcal H^{n-1}(\xi)\,dr\\[10pt]
&\le C_n\varepsilon_m^2
\int_1^R r^{n-1}
\sum_{j=1}^{N_m}
\int_{U_{m,j}}
\left|
\nabla_{\mathbb S^{n-1}}
\bigl(u(r\,\cdot)\bigr)(\xi)
\right|^2
\,d\mathcal H^{n-1}(\xi)\,dr\\[10pt]
&\le C_nR^2\varepsilon_m^2
\int_{\Omega_m\setminus B_1}|\nabla u|^2\,dx.
\end{aligned}
\]
Consequently, by density, we have
\eqref{eq:discussion-sun-poincare}.\medskip

We can now conclude the proof verifying the second Mosco condition.
Let $u_m\in H_0^1(\Omega_m)$ such that
converges weakly to $u$ in $H_0^1(D)$. Then, the sequence $(u_m)$ is bounded
in $H_0^1(D)$, and 
\eqref{eq:discussion-sun-poincare} gives
\[
\begin{aligned}
\int_{D\setminus B_1}|u_m|^2\,dx
&=\int_{\Omega_m\setminus B_1}|u_m|^2\,dx &\le C(n,R)\varepsilon_m^2
\int_{D}\abs{\nabla u_m}^2
\longrightarrow0.
\end{aligned}
\]
Thus $u_m$ converges strongly to zero in
$L^2(D\setminus B_1)$. On the other hand,
the weak convergence of $u_m$ to $u$ in $H_0^1(D)$ implies
weak convergence to $u$ in $L^2(D\setminus B_1)$.
By uniqueness of the weak limit, we conclude
that $u=0$ almost everywhere in $D\setminus B_1$, thus the second Mosco condition is therefore verified.\medskip

By \cite[Proposition~3.5.5]{HP18shape}, for equibounded open sets, the convergence in the sense of Mosco of $H^1_0(\Omega_m)$ to $H^1_0(B_1)$,
implies the strong convergence, in $H^1_0(D)$, of the torsion functions $w_{\Omega_m}$ to $w_{B_1}$. In particular, the Cauchy--Schwarz inequality yields
\[
\begin{aligned}
d_\gamma(\Omega_m,B_1)
&=\norm{w_{\Omega_m}-w_{B_1}}_{1}=\norm{w_{\Omega_m}-w_{B_1}}_{L^1(D)}\le |D|^{1/2}
\norm{w_{\Omega_m}-w_{B_1}}_{L^2(D)}
\longrightarrow0.
\end{aligned}
\]

Finally, by \autoref{thm-KJ-gamma}, we also have 
\[
\lim_{m}\delta_{\mathrm{KJ}}(\Omega_m)=0.\qedhere
\]\end{proof}

\paragraph{Funding}
The authors are members of GNAMPA and they were partially supported by Gruppo Nazionale per l’Analisi Matematica, la Probabilità e le loro Applicazioni (GNAMPA) of Istituto Nazionale di Alta Matematica (INdAM).

\printbibliography[heading=bibintoc]

\Addresses
\end{document}